\documentclass[12pt]{article}
\usepackage{latexsym}
\usepackage{amsmath}
\usepackage{amssymb}
\usepackage{amsthm}
\usepackage{amscd}
\usepackage[dvipdfmx]{graphicx}
\usepackage{color}
\input cyracc.def
\newfont{\tencyr}{wncyr10}

\begin{document}

\begin{center}
 {\large \bf  On mod $3$  triple Milnor invariants and triple cubic residue symbols in the Eisenstein number field}
 \end{center}

 \vspace{.02cm}
 
\begin{center}
Fumiya Amano, Yasushi Mizusawa and Masanori Morishita
\end{center}

\begin{center}
{\em Dedicated to Professor Takayuki Oda}
\end{center}

\vspace{.2cm}

{\small
 {\bf Abstract.} We introduce mod $3$ triple Milnor invariants and triple cubic residue symbols for certain primes of the Eisenstein number field $\mathbb{Q}(\sqrt{-3})$, following the analogies between knots and primes. Our triple symbol generalizes both  the cubic residue symbol and  R\'{e}dei's triple symbol, and describes the decomposition law of a prime in a mod $3$ Heisenberg extension of degree $27$ over $\mathbb{Q}(\sqrt{-3})$ with restricted ramification, which we construct concretely in the form similar to R\'{e}dei's dihedral extension over $\mathbb{Q}$.  We also give a cohomological interpretation of our symbols by triple Massey products in Galois cohomology. }

 \footnote[0]{
$\;\;$ Key words: Eisenstein number field, maximal pro-$3$ Galois group with restricted  \\
$\;\;$ ramification, mod 3 triple Milnor invariant,  mod 3 Heisenberg extension, triple cubic \\
$\;\;$ residue symbol,  triple Massey product\\
$\;\;$ The second author is partly supported by JSPS KAKENHI Grant Number JP26800010, \\
$\;\;$ Grant-in-Aid for Young Scientists (B)\\
$\;\;$ The third author is partly supported by JSPS KAKENHI Grant Number JP17H02837, \\
$\;\;$ Grant-in-Aid for Scientific Research (B)}

\vspace{.2cm}

\begin{center}
{\bf Introduction}
\end{center}

In this paper, we investigate a triple generalization of the cubic residue symbol in the Eisenstein number field $\mathbb{Q}(\sqrt{-3})$. It is L. R\'{e}dei who firstly studied such a generalization of the Legendre symbol in 1939 ([Rd]), aiming to generalize the  arithmetic of quadratic fields such as the theory of genera initiated by Gauss ([G]).  For distinct rational primes $p_1,p_2$ and $p_3$ satisfying $p_i \equiv 1 \mod 4$ and $\left(\frac{p_i}{p_j} \right) = 1 \; (1\leq i\neq j \leq 3)$, R\'{e}dei introduced a triple symbol $[p_1,p_2,p_3]$, which describes the decomposition law of $p_3$ in a dihedral extension $\frak{R}$, determined by $p_1$ and $p_2$, of degree $8$ over the rational number field  $\mathbb{Q}$.  Here  R\'{e}dei's extension $\frak{R}$ is given concretely by
$$ \frak{R} = \mathbb{Q}(\sqrt{p_1}, \sqrt{p_2}, \sqrt{\alpha}),  \leqno{(0.1)} $$
where $\alpha = x + y \sqrt{p_1}$ and $x, y$ are certain integers satisfying $x^2 - p_1y^2 - p_2 z^2 = 0$ with some non-zero integer $z$ (cf. [Rd]). It is characterized as the unique  Galois extension over $\mathbb{Q}$ whose Galois group is the dihedral group $D_8$ of order $8$ and which is unramified outside $p_1, p_2$ and the infinite prime with ramification index of each $p_i$ being $2$ ([A1]). It might not be clear, however, why such a dihedral extension and triple symbol should be considered as a natural generalization of a quadratic field and the Legendre symbol, and it seemed that his work had been overlooked for a long time (except some related works [F\"{o}], [Fu] etc). \\

In the late 1990s, M. Kapranov and the third author independently  interpreted the R\'{e}dei symbol as a mod $2$  arithmetic analogue of a triple linking number of a link, and further the third author introduced  mod $2$ arithmetic analogues for rational primes of the Milnor invariants (higher order linking numbers)  in link theory ([Mi2]), based on the analogies between primes and knots in {\em arithmetic topology} ([Ka], [Mo1]$\sim$[Mo4], [Rz]). For example, the mod $2$  Milnor invariant $\mu_2(12\cdots n) \in \mathbb{F}_2$ of length $n \geq 2$ for certain rational primes $p_1, p_2, \dots , p_n$ describes the decomposition law of $p_n$ in a certain nilpotent extension $\frak{K}(n)$, determined by $p_1, \dots , p_{n-1}$,  of degree $2^{\frac{1}{2}n(n-1)}$ over $\mathbb{Q}$, where the extension $\frak{K}(n)/\mathbb{Q}$ is unramified outside $p_1, p_2, \dots , p_{n-1}$ and the infinite prime with ramification index for each $p_i$ being $2$, and the Galois group ${\rm Gal}(\frak{K}(n)/\mathbb{Q})$ is isomorphic to the group $H_n(\mathbb{F}_2)$ consisting of $n$ by $n$ upper-triangular unipotent matrices  over $\mathbb{F}_2$.  For the cases that  $n = 2$ and $3$, where $H_2(\mathbb{F}_2) = \mathbb{Z}/2\mathbb{Z}$ and $H_3(\mathbb{F}_2) = D_8$, the mod $2$ Milnor invariants $\mu_2(12)$ and $\mu_2(123)$ give the Legendre symbol $\left( \frac{p_1}{p_2} \right)$ and the R\'{e}dei symbol $[p_1,p_2,p_3]$, respectively,  in the relations  
$$ (-1)^{\mu_2(12)} = \left( \frac{p_1}{p_2} \right), \;\;\;\; (-1)^{\mu_2(123)} = [p_1,p_2,p_3], \leqno{(0.2)}$$
and further it is shown  that  
$$\frak{K}(2) = \mathbb{Q}(\sqrt{p_1}), \;\;\;\; \frak{K}(3) = \frak{R}  \leqno{(0.3)}$$
([A1], [M1]$\sim$[M3]).  The first author introduced further the $4$-tuple quadratic residue symbol $[p_1,p_2,p_3,p_4] = (-1)^{\mu_2(1234)}$ by constructing concretely a candidate of $\frak{K}(4)$ ([A2]).  This unified interpretation may tell us that R\'{e}dei's dihedral extension and triple symbol would be a natural generalization of a quadratic field and the Legendre symbol.

However, these results concerns only rational primes, because they are based on an analogy between the structure of the maximal pro-$2$ Galois group over $\mathbb{Q}$ with restricted ramification and that of the group of a link in the $3$-sphere and there are obstructions to extend this analogy for number fields. We note that topological Milnor invariants are well defined for a link in an $l$-homology $3$-sphere ([Tu]), $l$ being a prime number, and that $l$-homology $3$-spheres correspond to number fields having trivial $l$-class group in the analogy of arithmetic topology. So it would be an interesting and important problem to overcome the difficulty concerning the obstructions and to introduce mod $l$ arithmetic Milnor invariants and multiple $l$-th power residue symbols for certain primes in such a number field containing a primitive $l$-th root of unity. 

In this paper, we introduce and investigate mod $3$ arithmetic triple Milnor invariants and triple cubic residue symbols for certain primes in the Eisenstein number field, by overcoming the difficulty concerning the obstruction mentioned above. In the following, let us describe our main results.
\\

Let $k$ be the Eisenstein number field $\mathbb{Q}(\zeta_3) = \mathbb{Q}(\sqrt{-3})$, $\zeta_3 := \frac{-1+\sqrt{-3}}{2}$. In this paper,  we shall introduce the mod $3$ triple Milnor invariant $\mu_3(123)$ for a set $S_0 = \{ \frak{p}_1,\frak{p}_2,\frak{p}_3\}$ of primes of $k$ satisfying ${\rm N}\frak{p}_i \equiv 1 \mod 9$  and $\left(\frac{\pi_i}{\pi_j}\right)_3 = 1 \; (1\leq i \neq j \leq 3)$ $(1\leq i \neq j  \leq 3)$, where $\pi_i$ is  the  unique prime element of $\frak{p}_i$ such that  $\pi_i \equiv 1 \; \mbox{mod}\; (3\sqrt{-3})$.  Although there is the obstruction $B_{S_0}$ for the Galois group $\frak{G}_{S_0}$ of the maximal pro-$3$ extension of $k$ unramified outside $S_0$ as mentioned above, we observe that if we add to $S_0$ a prime $\frak{p}_0$ with ${\rm N}\frak{p}_0 \equiv 4 \; \mbox{or}\; 7 \; \mbox{mod}\; 9$, the obstruction 
$B_S$ for $S = S_0\cup \{ \frak{p}_0\}$ vanishes so that the larger Galois group $\frak{G}_S$ of the maximal pro-$3$ extension $k_{S}$ of $k$ unramified outside $S$ has a  minimal presentation similar to a link group
$$ \frak{G}_{S} = \langle x_1, x_2, x_3 \; | \; x_1^{{\rm N}\frak{p}_1-1}[x_1,y_1] =  x_2^{{\rm N}\frak{p}_2-1}[x_2,y_2] = x_3^{{\rm N}\frak{p}_3-1}[x_3,y_3]=1  \rangle, \leqno{(0.4)} $$
where the letter $x_i$ represents a monodromy over $\frak{p}_i$ in $k_S/k$ and $y_i$ is a pro-$3$ word representing a Frobenius automorphism over $\frak{p}_i$ in $k_S/k$. Following the procedure in link theory ([Mi2]), we derive the mod $3$ Milnor numbers $\mu_3(ij)$ $(1\leq i \neq j \leq 3)$ and $\mu(123)$ from the presentation (0.4) by using the pro-$3$ Magnus expansion of $y_j$ and $y_3$, respectively. Then we firstly prove that $\mu_3(ij)$ is independent of the choice of $\frak{p}_0$ and so an invariant of the ordered pair ($\frak{p}_i$, $\frak{p}_j$), by noting that there exists (uniquely) a cyclic extension of degree $3$ over $k$ in which only $\frak{p}_i$ is ramified. As for $\mu_3(123)$, we shall prove, under $\mu_3(ij) = 0$ ($1\leq i \neq j \leq 3$), that 
$\mu_3(123)$  is independent of the choice of $\frak{p}_0$ and so  an invariant of the ordered triple ($\frak{p}_1, \frak{p}_2, \frak{p}_3$) if there exists (uniquely)  a Galois extension $\frak{K}_{\{ \frak{p}_1, \frak{p}_2\}}$ over $k$ whose Galois group is isomorphic to the mod $3$ Heisenberg group $H_3(\mathbb{F}_3)$ and in which only primes $\frak{p}_1, \frak{p}_2$  are ramified with ramification index of each $\frak{p}_i$ being $3$. We show that such a extension of $k$  exists if and only if  the class number of $k(\sqrt[3]{\pi_1\pi_2})$ is divisible by 9 and that the latter condition is satisfied if $\frak{p}_1$ and $\frak{p}_2$ are generated by  prime numbers. 

We then define the triple cubic symbol $[\frak{p}_1, \frak{p}_2,\frak{p}_3]_3$ by 
$$ [\frak{p}_1, \frak{p}_2,\frak{p}_3]_3 := \zeta_3^{\mu_3(123)}$$
as a cubic generalization of (0.2). In order to describe $[\frak{p}_1, \frak{p}_2,\frak{p}_3]_3$ arithmetically, we construct concretely the extension $\frak{K}_{\{\frak{p}_1,\frak{p}_2\}}$ in the form 
$$\frak{K}_{\{\frak{p}_1,\frak{p}_2\}} = k(\sqrt[3]{\pi_1}, \sqrt[3]{\pi_2}, \sqrt[3]{\theta}), \leqno{(0.5)}$$
under a certain assumption on $\theta$, where $\theta = x + y \sqrt[3]{\pi_1} + z(\sqrt[3]{\pi_1})^2 $ and $x, y, z$ are certain algebraic integers in $\mathbb{Z}[\zeta_3]$ satisfying $x^3 + \pi_1y^3 + \pi_1^2 z^3 - 3 \pi_1 xyz - \pi_2^3 w^3 = 0$ with some $w \in \mathbb{Z}[\zeta_3]$. The assumption on $\theta$ is proved to hold if $\frak{p}_1$ and $\frak{p}_2$ are generated by prime numbers. This extension  (0.5) was firstly constructed by the first author in [A3] and may be regarded as a cubic generalization of R\'{e}dei's extension $\frak{R}$ over $\mathbb{Q}$ in (0.1).  We can define the Artin symbol $\left( \frac{ \frak{K}_{\{\frak{p}_1,\frak{p}_2\}} /k}{\frak{p}_3}   \right)$ by the assumption and then we can show the formula
$$[\frak{p}_1, \frak{p}_2, \frak{p}_3]_3 = \frac{  \left(  \frac{ \frak{K}_{\{\frak{p}_1,\frak{p}_2\}}/k}{\frak{p}_3}   \right)(\sqrt[3]{\theta} )}{\sqrt[3]{\theta}}.$$
 So $[\frak{p}_1, \frak{p}_2,\frak{p}_3]_3$ describes the decomposition law of $\frak{p}_3$ in the extension $\frak{K}_{\{\frak{p}_1,\frak{p}_2\}} /k$. 

The power residue symbol is known to be described by the cup product in Galois cohomology ([Ko; 8.11], [Se; XIV, $\S 2$]). We generalize this relation for our triple cubic residue symbol by describing them by the Massey product in Galois cohomology. It is also an extension of the earlier works [Mo3], [V] in the case of the rational number field to the Eisenstein number field. 
 
Thus our results in this paper may form a first and natural generalization of the arithmetic of R\'{e}dei's triple symbols and cubic residue symbols, and also our work may have some importance in the respect that it would suggest a general theory of multiple power residue symbols in number fields, to be explored in the future.\\

Here are the contents of this paper.  In Section 1, we recall a theorem of Koch for the particular case of the maximal pro-$3$ Galois group $\frak{G}_S$ over the Eisenstein number field $k = \mathbb{Q}(\sqrt{-3})$ with given ramification primes $S$. We determine the obstruction $B_S$ and,  consequently, when $S$ contains a prime $\frak{p}$ with ${\rm N}\frak{p} \equiv 4 \; \mbox{or}\; 7 \; \mbox{mod}\; 9$, we describe a minimal presentation of $\frak{G}_S$ in the form similar to a link group.  In Section 2, for a given set of primes $S_0 = \{ \frak{p}_1, \dots , \frak{p}_r \}$ of $k$ with ${\rm N}\frak{p}_i \equiv 1 \; \mbox{mod}\; 9$, we derive the mod $3$ Milnor numbers $\mu_3(i_1\cdots i_n)$ $(1\leq i_1,\dots , i_n \leq r)$ from the Galois group $\frak{G}_S$ for $S = S_0 \cup \{ \frak{p}_0\}$ with ${\rm N}\frak{p}_0 \equiv 4 \; \mbox{or}\; 7 \; \mbox{mod}\; 9$ and show that they are invariants of the group $\frak{G}_S$. We introduce the notion of a R\'{e}dei type $H_n(\mathbb{F}_3)$-extension $\frak{K}_{\{ \frak{p}_1,\dots , \frak{p}_{n-1}\}}$ of $k$ for $\{ \frak{p}_1,\dots , \frak{p}_{n-1}\}$ with ${\rm N}\frak{p}_i \equiv 1 \; \mbox{mod}\; 9$. 
 In Section 3, we prove that the mod $3$ Milnor number $\mu_3(12)$ is independent of the choice of $\frak{p}_0$ and so an invariant of the ordered pair $(\frak{p}_1, \frak{p}_2)$, by noting the existence of the R\'{e}dei type $H_2(\mathbb{F}_3)$-extension of $k$ for $\{ \frak{p}_1\}$. We then show its relation with the cubic residue symbol. In Section 4,  we give a condition for a R\'{e}dei type $H_3(\mathbb{F}_3)$-extension $\frak{K}_{\{ \frak{p}_1, \frak{p}_2\}}$ of $k$ for $\{\frak{p}_1, \frak{p}_2\}$ to exist, and  show that the condition is satisfied if $\frak{p}_1$ and $\frak{p}_2$ are generated by prime numbers. We then prove that the mod $3$ triple Milnor number $\mu_3(123)$ is independent of the choice of $\frak{p}_0$ and so an invariant of the ordered triple $(\frak{p}_1, \frak{p}_2, \frak{p}_3)$ of primes in $S_k^{1 \, \mbox{{\scriptsize mod}}\, 9}$, under $\mu_3(ij) = 0$ ($1\leq i \neq j \leq 3$), if there exists (uniquely) a R\'{e}dei type $H_3(\mathbb{F}_3)$-extension $\frak{K}_{\{ \frak{p}_1, \frak{p}_2\}}$ of $k$ for $\{\frak{p}_1, \frak{p}_2\}$.  In Section 5, we construct concretely the extension  $\frak{K}_{\{ \frak{p}_1,  \frak{p}_{2} \}}$ in the form analogous to Redei's dihedral extension over $\mathbb{Q}$, under a certain assumption. We show that this assumption holds if $\frak{p}_1$ and $\frak{p}_2$ are generated by prime numbers. In Section 6, we introduce the triple cubic residue symbol in terms of the mod $3$ triple Milnor invariant, and describe it by using the concrete construction of $\frak{K}_{\{ \frak{p}_1,  \frak{p}_{2} \}}$.  In Section 7, we give a cohomological interpretation of our triple residue symbols in terms of  Massey products in Galois cohomology.\\
\\
{\it Notation.} For a number field $F$, ${\cal O}_F$ denotes the ring of integers of $F$. For a non-zero ideal $\frak{a}$ of ${\cal O}_F$, ${\rm N}\frak{a}$ denotes the norm of $\frak{a}$, ${\rm N}\frak{a} := \#({\cal   O}_F/\frak{a})$. We denote by $Cl(F)$ the ideal class group of $F$ and by $h_F$ the class number of $F$, $h_F := \#Cl(F)$.\\
For elements $a, b$ in a group $G$, the commutator $[a,b]$ is defined by $aba^{-1}b^{-1}$.  \\
Throughout this paper, let $\zeta_{3} := \frac{-1+\sqrt{-3}}{2}$ and let $k$ denote the Eisenstein number field, $k := \mathbb{Q}(\zeta_3) = \mathbb{Q}(\sqrt{-3})$. \\
We denote by $S_k^{\mbox{{\scriptsize non-}}3 }$ the set of finite primes of $k$ (maximal ideals of ${\cal O}_k$) which are not lying over $3$, and we set
$$ \begin{array}{l} S_k^{1\, \mbox{{\scriptsize mod}}\, 9 } := \{ \frak{p} \in S_k^{\mbox{{\scriptsize non-}}3 }\, | \, {\rm N}\frak{p} \equiv 1 \; \mbox{mod} \; 9 \},\\
S_k^{4, 7\, \mbox{{\scriptsize mod}}\, 9 } := \{ \frak{p} \in S_k^{\mbox{{\scriptsize non-}}3 }\, | \, {\rm N}\frak{p} \equiv 4 \; \mbox{or}\; 7 \; \mbox{mod} \; 9 \}.
\end{array}$$
\vspace{0.02cm}

\begin{center}
{\bf 1. Maximal pro-$3$ Galois groups with restricted ramification} 
\end{center}

In this section, we recall a theorem of Koch ([Ko]) on the maximal pro-$l$ Galois group $\frak{G}_S$ over a number field $k$ with given ramification primes $S$, for the particular  case where $l=3$ and $k$ is the Eisenstein number field $\mathbb{Q}(\sqrt{-3})$. For this case, we can determine the obstruction $B_S$ involved in a theorem of Koch and hence the generator rank of $\frak{G}_S$. Consequently, when $S$ contains a prime $\frak{p}$ with ${\rm N}\frak{p} \equiv 4 \; \mbox{or}\; 7 \; \mbox{mod}\; 9$, we can describe a presentation of $\frak{G}_S$ in the form similar to a link group ([Mi1], [Mi2]).\\

First, we note that ${\cal O}_k = \mathbb{Z}[\zeta_3]$, the unit group is ${\cal O}_k^{\times} = \{ \pm \zeta_3^e \, | \, e = 0, 1, 2\}$ and the class number $h_k$ is one.  Recall our notaion: $S_k^{\mbox{{\scriptsize non-}}3 }$ denotes the set of finite primes of $k$ which are not lying over $3$. For $\frak{p} \in S_k^{\mbox{{\scriptsize non-}}3 }$, we have ${\rm N}\frak{p} \equiv 1 \; \mbox{mod}\; 3$. We denote by  $S_k^{1\, \mbox{{\scriptsize mod}}\, 9 }$  (resp. $S_k^{4, 7\, \mbox{{\scriptsize mod}}\, 9 }$) the set of $\frak{p} \in S_k^{\mbox{{\scriptsize non-}}3 }$ satisfying ${\rm N}\frak{p} \equiv 1 \; \mbox{mod} \; 9$ (resp.  ${\rm N}\frak{p} \equiv 4 \; \mbox{or}\; 7 \; \mbox{mod} \; 9$), so that $S_k^{\mbox{{\scriptsize non-}}3 } = S_k^{1\, \mbox{{\scriptsize mod}}\, 9 } \sqcup S_k^{4, 7\, \mbox{{\scriptsize mod}}\, 9 }$. We start to show the following elementary lemma, which asserts that there is a standard prime element in each prime in  $S_k^{\mbox{{\scriptsize non-}}3 }$. It may be interesting, from the viewpoint of arithmetic topology, to remark that choosing a prime element in a finite prime is analogous to choosing an orientation of a knot. 
\\
\\
{\bf Lemma 1.1.}  {\em Let $\frak{p} \in S_k^{\mbox{{\scriptsize non-}}3 }$. Then $\frak{p} \in S_k^{1\, \mbox{{\scriptsize mod}}\, 9}$ if and only if there exists uniquely $\pi \in {\cal O}_k$ such that}
$$ \frak{p} = (\pi), \; \; \pi \equiv 1 \; \mbox{mod} \; (3\sqrt{-3}).$$
{\em Similarly, ${\rm N}\frak{p} \equiv 4$ $($resp. ${\rm N}\frak{p} \equiv 7) \; {\rm mod} \; 9$ if and only if there exists the unique prime element $\pi$ in  $\frak{p}$ satisfying $\pi \equiv 7 \; {\rm mod} \; (3\sqrt{-3})$ $($resp. $ \pi \equiv 4 \; {\rm mod} \; (3\sqrt{-3}))$. }\\
\\
{\em Proof.}  Suppose that $\frak{p} \in S_k^{1\, \mbox{{\scriptsize mod}}\, 9 }$. Since the class number of $k$ is one, there is $\pi' \in {\cal O}_k$ such that $\frak{p} = (\pi')$. Since ${\rm N}_{k/\mathbb{Q}}(\pi') > 0$, we have  
$${\rm N}_{k/\mathbb{Q}}(\pi') = {\rm N}\frak{p} \equiv 1 \mod 9.  \leqno{(1.1.1)}$$
Let $\frak{U} := ({\cal O}_k/(3\sqrt{-3}))^{\times}$. Let $\frak{l}$ be the unique maximal ideal of ${\cal O}_k$ lying over $3$, $\frak{l} = (1-\zeta_3)$. We note that $(3\sqrt{-3}) = \frak{l}^3$ and ${\cal O}_k/\frak{l} = \mathbb{Z}/3\mathbb{Z}$. So we have
$$ \begin{array}{ll} 
\frak{U}  & = \{ a_0 + a_1 \sqrt{-3} + a_2 (\sqrt{-3})^2 \; \mbox{mod}\; (3\sqrt{-3}) \; | \\
                                                   &   \;\;\;\;\;\;\;\;\;\;\;\;\;\;\;\; \;\;    \;\;\;\;\;\;\;\;\;\;\;\;\;\;\;\; \;\;   \;\;\;\;\;\;\;\;\; \;\;\;\;\; a_0 = 1,2,\; a_1, a_2 = 0, 1, 2 \}   \\
  & = \{ a + b \sqrt{-3} \; \mbox{mod}\; (3\sqrt{-3}) \; | \; a = 1,2,4,5,7,8, \; b = 0,1,2 \}.
\end{array} \leqno{(1.1.2)}
$$
Let $\frak{U}^{1}$ be the subgroup of $({\cal O}_k/(3\sqrt{-3}))^{\times}$ consisting of $\alpha$ mod $(3\sqrt{-3})$ with ${\rm N}_{k/\mathbb{Q}}(\alpha) \equiv 1 \; \mbox{mod} \; 9$. By the straightforward calculation using (1.1.2), we have
$$ \begin{array}{ll}
\frak{U}^1 & = \{ a + b \sqrt{-3} \; \mbox{mod}\; (3\sqrt{-3}) \; | \\
            &  \;\;\;\;\;\;\;\;\;\;\;\;\;\;\;\; \;\; (a, b) = (1,0), (8,0), (4,1), (4,2), (5,1), (5,2) \} \\
       &  = \langle -1 \; \mbox{mod}\; (3\sqrt{-3}) \rangle \times \langle \zeta_3 \; \mbox{mod}\; (3\sqrt{-3}) \rangle\\
       & = \{ \varepsilon \, \mbox{mod}\,  (3\sqrt{-3}) \; | \; \varepsilon \in {\cal O}_k^{\times} \},
            \end{array}
            \leqno{(1.1.3)}
$$
where $-1 \equiv 8$, $\zeta_3 \equiv 4+2\sqrt{-3}$ mod $(3\sqrt{-3})$.
By (1.1.1) and (1.1.3), there is a unit $\varepsilon \in {\cal O}_k^{\times}$ such that $\pi := \varepsilon \pi' \equiv 1$ mod $(3\sqrt{-3})$ and $(\pi) = (\pi')$. 

Suppose that $\frak{p} = (\pi) = (\varpi)$ and $\pi \equiv \varpi \equiv 1$ mod $(3\sqrt{-3})$. We can write  $\varpi = \varepsilon \pi$ for some $\varepsilon \in {\cal O}_k^{\times}.$  So $\varepsilon \equiv \varepsilon \pi = \varpi \equiv 1$ mod $(3\sqrt{-3})$. By (1.1.3), we must have $\varepsilon = 1$ and hence $\varpi = \pi$. 

Conversely, suppose that $ \frak{p} = (\pi), \; \pi \equiv 1 \mod (3\sqrt{-3})$. Writing $\pi = 1 + 3\sqrt{-3}(x+y\zeta_3)$ with $x, y \in \mathbb{Z}$, we see easily ${\rm N}\frak{p} = {\rm N}_{k/\mathbb{Q}}(\pi) \equiv 1 \mod 9$.

Similarly, the latter assertions are verified by a straightforward computation using (1.1.2). $\;\; \Box$
\\

Let $S$  be a finite set of $s$ distinct  primes $\frak{p}_1,\dots , \frak{p}_s$ in $S_k^{\mbox{{\scriptsize non-}}3 }$.  Let $k_S$ denote the maximal pro-$3$  extension  of  $k$, unramified outside $S$,  in a fixed algebraic closure $\bar{k}$. Let $\frak{G}_{S}$ denote the Galois group of $k_S$ over $k$, $\frak{G}_{S} := {\rm Gal}(k_S/k)$. We describe the structure of the pro-$3$ group $\frak{G}_{S}$ in a certain unobstructed case. For this, we first recall a result due to  Iwasawa on local Galois groups ([Iw2]). For each $i$ ($1\leq i \leq s$), let $k_{\frak{p}_i}$ be the $\frak{p}_i$-adic field which has the unique prime element $\pi_i$ as given in Lemma 1.1. We fix an algebraic closure $\overline{k}_{\frak{p}_i}$ of $k_{\frak{p}_i}$ and an embedding $\bar{k} \hookrightarrow \overline{k}_{\frak{p}_i}$. Let $\tilde{k}_{\frak{p}_i}$ denote the maximal pro-$3$ extension of $k_{\frak{p}_i}$ in $\overline{k}_{\frak{p}_i}$ and $\frak{G}_{\frak{p}_i}$ denote the Galois group of $\tilde{k}_{\frak{p}_i}$ over $k_{\frak{p}_i}$, $\frak{G}_{\frak{p}_i} := {\rm Gal}(\tilde{k}_{\frak{p}_i}/k_{\frak{p}_i})$. Then we have
$$ \tilde{k}_{\frak{p}_i} = k_{\frak{p}_i}( \zeta_{3^n}, \sqrt[3^n]{\pi_i} \; | \; n \geq 1 ),$$
where $\zeta_{3^n}$ denotes a primitive $3^n$-th root of unity in $\bar{k}$  such that $\zeta_3 = \frac{-1+\sqrt{-3}}{2}$ and $(\zeta_{3^a})^{3^b} = \zeta_{3^{a-b}}$ for all $a \geq b$.  The local Galois group $\frak{G}_{\frak{p}_i}$  is  then topologically generated by the monodromy $\tau_i$ and (an extension of) the Frobenius automorphism $\sigma_i$ which are defined by
$$ \begin{array}{ll}
\tau_i(\zeta_{3^n}) := \zeta_{3^n}, & \tau_i(\sqrt[3^n]{\pi_i}) := \zeta_{3^n}\sqrt[3^n]{\pi_i},\\
\sigma_i(\zeta_{3^n}) := \zeta_{3^n}^{{\rm N}\frak{p}_i},& \sigma_i(\sqrt[3^n]{\pi_i}) := \sqrt[3^n]{\pi_i}
\end{array}
\leqno{(1.2)}
$$
and subject to the relation
$$ \tau_i^{{\rm N}\frak{p}_{i} -1}[\tau_i, \sigma_i] = 1. \leqno{(1.3)}$$
  For each $i$ ($1\leq i \leq s$), the fixed  embedding $\bar{k} \hookrightarrow \overline{k}_{\frak{p}_i}$ gives an embedding $k_S \hookrightarrow \tilde{k}_{\frak{p}_i}$, hence a prime $\frak{P}_i$ of $k_S$ lying over $\frak{p}_i$. 
We denote by the same letters $\tau_i$ and $\sigma_i$ the images of $\tau_i$ and $\sigma_i$, respectively, under the homomorphism 
$$\varphi_{\frak{p}_i, S} : \frak{G}_{\frak{p}_i} \longrightarrow \frak{G}_{S}   \leqno{(1.4)}$$ 
induced by the embedding $k_S \hookrightarrow \tilde{k}_{\frak{p}_i}$. Then $\tau_i$ is a topological generator of the inertia group of the prime $\frak{P}_i$ and $\sigma_i$ is an extension of the Frobenius automorphism of the maximal subextension of $k_S/k$ for which $\frak{P}_i$ is unramified. We call simply $\tau_i$ and $\sigma_i$  a {\it monodromy over $\frak{p}_i$} in $k_S/k$ and  a {\it Frobenius automorphism over $\frak{p}_i$} in $k_S/k$, respectively. We note that the restriction of $\tau_i$ to the maximal abelian subextension $k_S^{\rm ab}$ of $k_S/k$ is given by the Artin symbol
$$ \tau_i|_{k_S^{\rm ab}} = ( \tilde{\xi}_i, k_S^{\rm ab}/k) \leqno{(1.5)}$$
for an idele $\tilde{\xi}_i$ of $k$ whose $\frak{p}_i$-component is a primitive $({\rm N}\frak{p}_i - 1)$-th root of unity $\xi_i$ in $k_{\frak{p}_i}^{\times}$ and other components are all 1.

Since the ideal class group $Cl(k)$  is trivial, class field theory tells us that the monodromies $\tau_1, \dots , \tau_s$ generate topologically the global Galois group $\frak{G}_{S}$. However, they may not be a minimal set of generators  in general. In fact, noting that $k$ contains $\zeta_3$, Shafarevich's theorem ([Ko; Satz 11.8]) tells us that the minimal number $d(\frak{G}_{S})$ of generators of $\frak{G}_{S}$ is given by
 $$ d(\frak{G}_{S}) = s - 1 + \dim_{\mathbb{F}_3} B_S.  \leqno{(1.6)}$$
Here the obstruction $B_S$ is defined by 
$$B_S := \{ a \in k^{\times} \; | \; (a) = \frak{a}^3,  \;  a \in (k_{\frak{p}_i}^{\times})^3  \; \mbox{for all}\; 1 \leq i \leq s \}/(k^{\times})^3,  \leqno{(1.7)}$$
where $\frak{a}$ is a fractional ideal of ${\cal O}_k$. For the Eisenstein number field, we can determine $B_S$ and hence $d_S(\frak{G}_S)$ as follows.\\
\\
{\bf Proposition 1.8.} {\it If $S$ is a subset of $S_k^{1\, \mbox{{\scriptsize mod}}\, 9}$, then $B_S = \langle \zeta_3 \, {\rm mod}\, (k^{\times})^3 \rangle \simeq  \mathbb{F}_3$ and hence $d(\frak{G}_S)= s$, namely, $\tau_1,\dots , \tau_s$ are minimal generators of $\frak{G}_{S}$. If $S$ contains a prime $\frak{p}$ in $S_k^{4, 7\, \mbox{{\scriptsize mod}}\, 9}$, then $B_S = \{ 1\}$ and hence $d(\frak{G}_S) = s-1$, namely, one of $\tau_1,\dots , \tau_s$ is redundant for minimal generators of $\frak{G}_{S}$.}\\
\\
{\em Proof.}  Suppose that  $\frak{p}_i \in  S_k^{1\, \mbox{{\scriptsize mod}}\, 9 }$ for $1\leq i \leq s$. Let $a \in k^{\times}$ satisfy  $(a) = \frak{a}^3$ and $a \in (k_{\frak{p}_i}^{\times})^3$ for $1\leq i \leq s$. Writing $\frak{a} = (b)$ with $b \in k^{\times}$, we have $a = \varepsilon b^3$ with $\varepsilon \in {\cal O}_k^{\times} = \{ \pm \zeta_3^e \, | \, e = 0, 1, 2\}$ and so $a \equiv \zeta_3^e$ mod $(k^{\times})^3$ for some $e = 0,1,2$. Noting that $\zeta_3 \in (k_{\frak{p}_i}^{\times})^3$ for all $1\leq i \leq s$ and $\zeta_3 \notin (k^{\times})^3$, we see that $\zeta_3$ mod $(k^{\times})^3$ is a basis of $B_S$.
\\
Suppose that  there is $\frak{p}_j \in S_k^{4, 7\, \mbox{{\scriptsize mod}}\, 9 }$. Let $a \in k^{\times}$ satisfy  $(a) = \frak{a}^3,  \;  a \in (k_{\frak{p}_i}^{\times})^3$ for  $1\leq i \leq s$. As in the above, writing $\frak{a} = (b)$, we have $a \equiv \zeta_3^e$ mod $(k^{\times})^3$ for some $e = 0,1,2$.  Since $\frak{p}_j$ is inert in $\mathbb{Q}(\zeta_9)/k$, $\zeta_9 \notin k_{\frak{p}_j}^{\times}$ and so $\zeta_3 \notin (k_{\frak{p}_j}^{\times})^3$. Therefore $e = 0$ and so $a \in (k^{\times})^3$. Hence $B_S = \{ 1 \}$. 
\\
The assertion for minimal generators of $\frak{G}_S$ follows from (1.6) $\;\; \Box$\\
\\
As for a redundant generator, class field theory tells us the following refined information. \\
\\
{\bf Proposition 1.9.} {\it Assume that $S$ contains  a prime $\frak{p}_j$ in $S_k^{4, 7\, \mbox{{\scriptsize mod}}\, 9 }$. Then we can exclude the monodromy $\tau_j$ over $\frak{p}_j$ from $\tau_1, \dots , \tau_s$ to obtain minimal generators of  $\frak{G}_{S}$.}
\\
\\
{\it Proof.} We may suppose that $\frak{p}_1, \dots, \frak{p}_r \in S_k^{1\, \mbox{{\scriptsize mod}}\, 9 }$ and $\frak{p}_{r+1}, \dots, \frak{p}_s \in S_k^{4, 7\, \mbox{{\scriptsize mod}}\, 9}$.
For $g \in \frak{G}_{S}$, we let $[g] := g$ mod $\frak{G}_{S}^3[\frak{G}_S, \frak{G}_S]$. We need to show that we can exclude one $[\tau_j]$ with $r+1\leq j \leq s$ from $[\tau_1],\dots , [\tau_s]$ to obtain a basis of the Frattini quotient $\frak{G}_S/\frak{G}_{S}^3[\frak{G}_S, \frak{G}_S]$. 
For a prime $\frak{p}$ of $k$, let ${\cal O}_{\frak{p}}$ denote the ring of $\frak{p}$-adic integers of $k_{\frak{p}}$. (For the infinite prime $\frak{p}_{\infty}$, we set ${\cal O}_{\frak{p}_{\infty}}:= k_{\frak{p}_{\infty}} = \mathbb{C}$). Let $J_k$ be the idele group of $k$. Let $U_k$ be the subgroup of $J_k$ consisting of unit ideles of $k$, $U_k := \prod_{\frak{p}} {\cal O}_{\frak{p}}^{\times}$, and let $U_{S}$ denote the subgroup of $U_k$ whose $\frak{p}$-component is $1$ for $\frak{p} \in S$. By class field theory, we have the canonical isomorphism by the Artin symbol of $\mathbb{F}_3$-vector spaces 
$$ J_k/U_{S}J_k^3 k^{\times} \; \simeq \; \frak{G}_{S}/\frak{G}_{S}^3[\frak{G}_S, \frak{G}_S].  \leqno{(1.9.1)}$$
Since the ideal class group $Cl(k) = J_k/U_kk^{\times}$ is trivial and $B_S = \{1\}$ by Proposition 1.8, we have the following exact sequence of $\mathbb{F}_3$-vector spaces (cf. [Ko; (11.11)])
$$  0 \longrightarrow {\cal O}_k^{\times}/({\cal O}_k^{\times})^3  \stackrel{\delta}{\longrightarrow} \prod_{i=1}^s {\cal O}_{\frak{p}_i}^{\times}/({\cal O}_{\frak{p}_i}^{\times})^3 \stackrel{\iota}{\longrightarrow} J_k/U_{S}J_k^3 k^{\times} \longrightarrow 0,  \leqno{(1.9.2)}$$
where $\delta$ is the diagonal map and $\iota$ is induced by the natural inclusion $\prod_{i=1}^s {\cal O}_{\frak{p}_i}^{\times} \hookrightarrow J_k$.
 Let $\xi_i$ be a primitive $({\rm N}\frak{p}_i -1)$-th root of unity in ${\cal O}_{\frak{p}_i}$ $(1\leq i\leq s)$. We denote by $[\xi_i]$ the element of $\prod_{i=1}^{s} {\cal O}_{\frak{p}_i}^{\times}/( {\cal O}_{\frak{p}_i}^{\times})^3$ whose $\frak{p}_i$-component is $\xi_i$ mod $( {\cal O}_{\frak{p}_i}^{\times})^3$ and other components are all 1, so that $[\xi_1], \dots , [\xi_s]$ form a basis of the $\mathbb{F}_3$-vector space $\prod_{i=1}^{s} {\cal O}_{\frak{p}_i}^{\times}/( {\cal O}_{\frak{p}_i}^{\times})^3$. We let $[\zeta_3] := \delta(\zeta_3 \; \mbox{mod}\; ({\cal O}_k^{\times})^3)$. Then we have an equation in $\prod_{i=1}^{s} {\cal O}_{\frak{p}_i}^{\times}/( {\cal O}_{\frak{p}_i}^{\times})^3$
$$ [\zeta_3] = \sum_{i=1}^s a_i [\xi_i] \;\; (a_i \in \mathbb{F}_3).  \leqno{(1.9.3)}$$
Note that the isomorphism (1.9.1) sends $\iota([\xi_i])$ to $[\tau_i] := \tau_i$ mod $\frak{G}_{S}^3[\frak{G}_S, \frak{G}_S]$ by (1.5) and sends $\iota([\zeta_3])$ to 0 by (1.9.2). Hence  we obtain from (1.9.3) the equation in $\frak{G}_{S}/\frak{G}_{S}^3[\frak{G}_S, \frak{G}_S]$
$$ 0 = \sum_{i=1}^s a_i [\tau_i]. \leqno{(1.9.4)}$$
\\
Let $1 \leq i \leq r$. Looking at the $\frak{p}_i$-component of (1.9.3), we have $\zeta_3 \equiv \xi_i^{a_i}$ mod $( {\cal O}_{\frak{p}_i}^{\times})^3$. Since ${\rm N}\frak{p}_i \equiv 1$ mod 9, we have $\zeta_3 \in ( {\cal O}_{\frak{p}_i}^{\times})^3$ and hence $a_i = 0$.\\
Let $r+1 \leq j \leq s$. Looking at the $\frak{p}_j$-component of (1.9.3), we have $\zeta_3 \equiv \xi_j^{a_j}$ mod $( {\cal O}_{\frak{p}_j}^{\times})^3$.  Since ${\rm N}\frak{p}_j \equiv 4$ or $7 \; \mbox{mod} \; 9$, we can write ${\rm N}\frak{p}_j -1 = 3m, (3,m)=1$. So we have  $\xi_j^{a_j m} = \zeta_3^m \neq 1$ and hence $a_j \neq 0$.  Therefore we have from (1.9.4)
$$ 0 = \sum_{j=r+1}^s a_i [\tau_j] \;\; (a_j \in \mathbb{F}_3^{\times}),$$
from which we can exclude arbitrary one $[\tau_j]$ with $r+1\leq j  \leq s$ from $[\tau_1],\dots , [\tau_s]$ to obtain a basis of $\frak{G}_S/\frak{G}_{S}^3[\frak{G}_S, \frak{G}_S]$.  $\;\; \Box$ \\
\\
Suppose that $B_S = \{1\}$. Noting that $k$ contains $\zeta_3$, it is shown by combining [Ko; Satz 6.11] ([Ko; Satz 6.14]), [Ko; Satz 11.3] and  [Ko; Satz 11.4] that the relations    of the Galois group $\frak{G}_S$ are given by the local relations (1.3) for minimal generators. Summing up, we have the following.\\
\\
{\bf Theorem 1.10} (Koch). {\it  Assume that $S$ contains a prime $\frak{p}_j$ in $S_k^{4, 7\, \mbox{{\scriptsize mod}}\, 9}$. Then the pro-$3$ group $\frak{G}_{S}$ has the following minimal presentation
$$\begin{array}{l} \frak{G}_{S} = \langle \,  x_{1}, \dots, x_{j-1}, x_{j+1}, \dots, x_s  \; | \; x_{1}^{{\rm N}\frak{p}_{1} -1}[x_{1},y_{1}] = \cdots = x_{j-1}^{{\rm N}\frak{p}_{j-1} -1}[x_{j-1},y_{j-1}] \\
\;\;\;\;\;\;\;\;\;\;  \;\;\;\;\; \; \;\;\;\;\; \;\;\;\;\; \;\;\;\;\; \;\;\;\;\;  \;\;\;\;\; = x_{j+1}^{{\rm N}\frak{p}_{j+1} -1}[x_{j+1},y_{j+1}] = \cdots =  x_{s}^{{\rm N}\frak{p}_{s} -1}[x_{s},y_{s}] =1 \, \rangle,
\end{array}$$
where  $x_{i}$ denotes the letter which represents a monodromy $\tau_{i}$ over $\frak{p}_{i}$ in $k_S/k$ and $y_{i}$ denotes the free pro-$3$ word of $x_1, \dots, x_{j-1}, x_{j+1}, \dots , x_{s}$  which represents 
a Frobenius automorphism over $\frak{p}_{i}$ in $k_S/k$.}\\

\begin{center}
{\bf 2. Mod 3  Milnor numbers}
\end{center}

In this section, for a given finite subset $S_0 = \{ \frak{p}_1, \dots , \frak{p}_r \}$ of $S_k^{1\, \mbox{{\scriptsize mod}}\, 9 }$, we introduce the mod $3$ Milnor numbers $\mu_3(i_1\cdots i_n)$ $(1\leq i_1,\dots , i_n \leq r)$ by using Theorem 1.10 applied to $S = S_0 \cup \{ \frak{p}_0\}$ with $\frak{p}_0 \in S_k^{4, 7\, \mbox{{\scriptsize mod}}\, 9}$, and we prove that the mod $3$ Milnor numbers are invariants of the pro-$3$ Galois group $\frak{G}_S$.\\

First, we recall the pro-$l$ Magnus expansion for a free pro-$l$ group $\frak{F}_N$ on letters $x_1,\dots , x_N$, where $l$ is a prime number (cf. [Ko; 7], [Mo4; 8.3]).  Let $\mathbb{F}_l[[\frak{F}_N]]$ denote the completed group algebra of $\frak{F}_N$ over $\mathbb{F}_l$ and let $\epsilon_{\mathbb{F}_l[[\frak{F}_N]]} : \mathbb{F}_l[[\frak{F}_N]] \rightarrow \mathbb{F}_l$ be the augmentation homomorphism with the augmentation ideal $I_{\mathbb{F}_l[[\frak{F}_N]]} := {\rm Ker}(\epsilon_{\mathbb{F}_l[[\frak{F}_N]]})$.  Let $\mathbb{F}_3 \langle \langle X_1,\dots , X_N \rangle \rangle$ denote the formal power series algebra over $\mathbb{F}_l$ in non-commuting variables $X_1, \dots , X_N$. Sending $x_i$ to $1 + X_i$ for $1 \leq i \leq r$, we have the (pro-$l$) {\em Magnus isomorphism} of topological $\mathbb{F}_l$-algebras
$$ \Theta_N : \mathbb{F}_l[[\frak{F}_{N}]]  \; \stackrel{\sim}{\longrightarrow} \;    \mathbb{F}_l \langle \langle X_1, \dots , X_{N} \rangle \rangle. $$
 For $\alpha \in \mathbb{F}_l[[\frak{F}_N]]$, $\Theta_N(\alpha)$ is called the {\it Magnus expansion} of $\alpha$ over $\mathbb{F}_l$. For a multi-index $I = (i_1 \cdots i_n)$,  $1\leq i_1,\dots , i_n \leq r$, we set $ |I| := n \; \mbox{and} \; X_I  := X_{i_1} \cdots X_{i_n}.$ 
We denote the coefficient of $X_I$ in the Magnus expansion $\Theta_N(\alpha)$ by  $\mu_l(I;\alpha)$, called the  mod $l$ {\em Magnus coefficient} of $\alpha$ for $I$, so that we have
$$ \Theta_N(\alpha) = \epsilon_{\mathbb{F}_l[[\frak{F}_N]]}(\alpha) + \sum_{|I|\geq 1} \mu_l(I;\alpha) X_I.$$
In terms of the pro-$l$ Fox free derivatives $\frac{\partial}{\partial x_i} : \mathbb{F}_l[[\frak{F}_N]] \rightarrow \mathbb{F}_l[[\frak{F}_N]]$ over $\mathbb{F}_l$ (cf. [Ih; \S 2], [Mo4; 8.3], [O]), we can write $\mu_l(I;\alpha)$ for $I = (i_1\cdots i_n)$ as
$$ \mu_l(I;\alpha) = \epsilon_{\mathbb{F}_l[[\frak{F}_N]]}\left(\frac{ \partial^n \alpha}{\partial x_{i_1} \cdots \partial x_{i_n}}\right).$$
Here are some basic properties of mod $l$ Magnus coefficients.\\
(2.1.1) For $\alpha, \beta \in \mathbb{F}_l[[\frak{F}_N]]$ and a multi-index $I$, we have
$$ \mu_l(I;\alpha \beta) = \sum_{I = JK} \mu_l(J;\alpha)\mu_l(K;\beta),$$
where the sum ranges over all pairs $(J,K)$ of multi-indices such that $JK = I$.\\
(2.1.2) ({\it Shuffle relation}) For $ f \in \frak{F}_N$ and multi-indices $I, J$ with $|I|, |J| \geq 1$, we have
$$ \mu_l(I;f)\mu_l(J;f) = \sum_{H \in {\rm Sh}(I,J)} \mu_l(H;f),$$
where ${\rm Sh}(I,J)$ denotes the set of the results of all shuffles of $I$ and $J$ ([CFL], [Mo4; 8.1]).\\
(2.1.3) For $f \in \frak{F}_N$ and $d \geq 2$, we have
$$ \mu_l(I;f) = 0 \; \mbox{for} \; |I| < d \; \Longleftrightarrow \; f \in \frak{F}_N^{(d)}, $$
where $\frak{F}_N^{(d)} := \frak{F}_M \cap (1+I_{\mathbb{F}_l[[\frak{F}_N]]}^d)$, the $d$-th term of the mod $l$ Zassenhaus filtration of $\frak{F}_N$. It is known that $\frak{F}_N^{(d)} = (\frak{F}_N^{([d/l])})^l \prod_{i+j=d}[\frak{F}_N^{(i)},\frak{F}_N^{(j)}]$, where $[d/l]$ stands for the Gauss symbol ([DDMS; Definition 11.1, Theorem 12.9]).\\

Let us be back in our arithmetic situation, where $l = 3$. Let $S_0 := \{ \frak{p}_1, \dots \, \frak{p}_r \}$ be a given set of $r$ distinct primes in $S_k^{1\, \mbox{{\scriptsize mod}}\, 9}$. Let us choose  a  prime $\frak{p}_0 \in S_k^{4, 7\, \mbox{{\scriptsize mod}}\, 9}$  and we let $S := S_0 \cup \{\frak{p}_0 \} = \{ \frak{p}_0, \frak{p}_1, \dots , \frak{p}_r\}$.   Let $\frak{G}_S$ denote the Galois group of the maximal pro-$3$ extension $k_S$ over $k$ unramified outside $S$. Let $x_i$ denote a letter representing a monodromy $\tau_i$ over $\frak{p}_i$ in $k_S/k$ $(1\leq i \leq r)$ and let $\frak{F}_r$ denote the free pro-$3$ group on $x_1, \dots , x_r$. Let $y_i$ denote the free pro-$3$ word which represents a Frobenius automorphism $\sigma_i$ over $\frak{p}_i$ in $k_S/k$. Let $\frak{N}_S$ denote the closed subgroup of $\frak{F}_r$ generated normally by $x_{1}^{{\rm N}\frak{p}_{1} -1}[x_{1},y_{1}], \dots ,  x_{r}^{{\rm N}\frak{p}_{r} -1}[x_{r},y_{r}]$. By Theorem 1.10, we have the minimal presentation
$$ \begin{array}{ll} \frak{G}_{S}  & = \langle \,  x_{1}, \dots, x_{r} \; | \; x_{1}^{{\rm N}\frak{p}_{1} -1}[x_{1},y_{1}] = \cdots =  x_{r}^{{\rm N}\frak{p}_{r} -1}[x_{r},y_{r}] =1 \rangle \\
& = \frak{F}_r/\frak{N}_S.  \end{array} \leqno{(2.2)}$$

 For a multi-index $I = (i_1\cdots i_n)$ $(1\leq i_1,\dots , i_n \leq r)$, we define the mod 3 {\it Milnor number} $\mu_3(I) = \mu_3^S(I)$ with respect to $S$ by the mod 3 Magnus coefficient of $y_{i_n}$ for $I' := (i_1\cdots i_{n-1})$ 
$$ \mu_3(I) := \mu_3(I'; y_{i_n}),  \leqno{(2.3)}$$
and we set $\mu_3(I) := 0$ if $|I| = 1$.  Let $e_I$ be the maximal integer $e$ satisfying ${\rm N}\frak{p}_i \equiv 1$ mod $3^e$ for all $i = i_1, \dots , i_n$ and we set $m_I := 3^{e_I}$. Note that $e_I \geq 2$, since ${\rm N}\frak{p}_i \equiv 1 \; \mbox{mod}\; 9$ for $1\leq i \leq r$. \\
\\
{\bf Theorem 2.4.} {\it  Let $I = (i_1\cdots i_n)$ be a multi-index satisfying $1\leq i_1,\dots ,i_n \leq r$ and $2 \leq |I|  \leq m_I$. Assume that $\mu_3(J) = 0$ if $|J| < |I|$. Then $\mu_3(J) = 0$ and $\mu_3(I)$ are independent of a choice of a monodromy $\tau_i$ and a Frobenius $\sigma_i$ over $\frak{p}_i$ for $1\leq i \leq r$, namely, a choice of a prime $\frak{P}_i$ in $k_S$ lying over $\frak{p}_i$ (equivalently an embedding $k_{\frak{p}_i}(l) \hookrightarrow k_{S}$), and so $\mu_3(I)$ is an invariant of the Galois group $\frak{G}_S$.} \\
\\
{\it Proof.} 
Note that $m_I \le m_J$. We see that by induction, $\mu_3(J) = 0$ is independent of a choice of $\tau_i$ and $\sigma_i$.
Since $\frak{G}_S$ has a presentation (2.2), we need to verify 
\vskip 2pt

(i) $\mu_3(I)$ is unchanged when $x_i$ is replaced by its conjugate for $1\leq i \leq r$.
\vskip 2pt

(ii) $\mu_3(I)$ is unchanged when $y_{i_n}$ is replaced by its conjugate in $\frak{F}_{r}$.
\vskip 2pt

(iii) $\mu_3(I)$ is unchanged when $y_{i_n}$ is multiplied by a product of conjugates \\
$\;\;\;\;\;\;\;\;\;\; $ of $(x_i^{{\rm N}\frak{p}_i - 1}[x_i,y_i])^e$ for $1\leq i \leq r,   e = \pm 1$.
\vspace{.04cm}\\
We set $I' := (i_1\cdots i_{n-1})$.

Proof of (i).  Suppose that $x_i$ is replaced by $x_i^{\dag} = x_jx_ix_j^{-1}$ ($1\leq i, j \leq r$). Since $x_i = x_j^{-1}x_i^{\dag}x_j$, we have $X_i = (1 - X_j + X_j^2 - \cdots)X_i^{\dag}(1 + X_j)$ and hence
$$ X_i = X_i^{\dag} + ( \mbox{terms involving} \; X_jX_i^{\dag} \; \mbox{or} \; X_i^{\dag} X_j).  \leqno{(2.4.1)} $$
Each time $X_i$  appears in the Magnus expansion $\Theta_{r}(y_{i_n}) = 1 + \sum_{J} \mu_3(J i_{n})X_J$, it is to be replaced by the above expansion (2.4.1) and we finally reach the new expansion of $\Theta_{r}(y_{i_n})$ in $X_1^{\dag}, \dots , X_{r}^{\dag}$, by which we denote $\Theta_{r}^{\dag}(y_{i_n})$. Then we can easily see that the coefficient of $X_{i_1}^{\dag} \cdots X_{i_{n-1}}^{\dag}$ in $\Theta_{r}^{\dag}(y_{i_n})$, denoted by $\mu_3^{\dag}(I)$, is of the form
$$ \mu_3^{\dag}(I) = \mu_3(I) + \sum_{J} \mu_3(Ji_{n}),$$
where $J$ runs over some proper subsequences of $I'$. Therefore, by the assumption, we have
$$ \mu_3^{\dag}(I)  = \mu_3(I).$$
Similarly, $\mu_3(I)$ is proved to be unchanged when $x_i$ is replaced by $x_j^{-1}x_i x_j$  ($1 \leq i, j \leq r$). So $\mu_3(I)$  is unchanged when $x_i$ is replaced by its conjugate in the (discrete) free group $F_{r}$ generated by $x_1,\dots , x_r$. Since $F_{r}$ is dense in $\frak{F}_{r}$ and $\mu_3(I)$ takes discrete values, $\mu_3(I)$ is unchanged when $x_i$ is replaced by its conjugate in $\frak{F}_{r}$. 

Proof of (ii). By comparing the coefficients of $X_{I'}$ in the both sides of the equality 
$$ \Theta_{r}(x_i y_{i_n}x_i^{-1}) = (1+X_i)\Theta_{r}(y_{i_n})(1 - X_i + X_i^2 - \cdots) $$
for $1 \leq i \leq r$ and by the assumption, we have
$$ \mu_3(I'; x_iy_{i_n} x_i^{-1}) = \mu_3(I). $$
Similarly, we have $\mu_3(I'; x_i^{-1} y_{i_n} x_i) = \mu_3(I).$  By the same argument as in the proof of (i), the assertion (ii) is proved. 

Proof of (iii).  Let $J$ be any subsequence of $I'$, $1\leq i \leq r$ and $e = \pm 1$. We will prove that
$$ \mu_3( J; (x_i^{{\rm N}\frak{p}_i -1} [x_i,y_i] )^e) = 0.  \leqno{(2.4.2)}$$
First we prove that 
$$ \mu_3(J; [x_i, y_i]^e) = 0. \leqno{(2.4.3)}$$
Comparing the coefficients of $X_J$ in the equality
$$ \Theta_{r}([x_i,y_i]^e) = \left\{ \begin{array}{l}
 1 + (\Theta_{r}(x_iy_i) - \Theta_{r}(y_ix_i))\Theta_{r}(x_i^{-1})\Theta_{r}(y_i^{-1}) \;\; (e = 1)\\
 1 + (\Theta_{r}(y_ix_i) - \Theta_{r}(x_iy_i))\Theta_{r}(y_i^{-1})\Theta_{r}(x_i^{-1}) \;\; (e = -1),
 \end{array}\right.
$$
we have
$$  \mu_3(J;[x_i,y_i]^e)  = e(\mu_3(J;x_iy_i) - \mu_3(J;y_ix_i))
                                                                 + \sum_A ( \mu_3(A;x_iy_i) - \mu_3(A;y_ix_i) )c_A,$$
where $A$ runs over some proper subsequences of $J$ and $c_A \in \mathbb{F}_3$. So, in order to prove (2.4.3), it is enough to show that for any subsequence $J$ of $I'$ and $1\leq i \leq r$,
$$ \mu_3(J;x_iy_i) - \mu_3(J; y_ix_i)  = 0.  \leqno{(2.4.4)}$$
Let $J = (j_1\cdots j_a)$. By the straightforward computaion, we have
$$ \mu_3(J; x_iy_i) = \left\{ \begin{array}{ll} 
\mu_3(Ji) \;\; &  (i \neq j_1),\\
\mu_3(Jj_1) + \mu_3(j_2\cdots j_a j_1) \;\; & (i = j_1),
\end{array} \right.
$$
and
$$
\mu_3(J; y_ix_i) = \left\{ \begin{array}{ll} 
\mu_3(Ji) \;\; & (i \neq j_a),\\
\mu_3(Jj_a) + \mu_3(J) \;\; & (i = j_a).
\end{array} \right.
$$
Hence we have
$$ \begin{array}{ll} \mu_3(J; x_iy_i) - \mu_3(J; y_ix_i) & = 
\left\{ \begin{array}{ll}
\mu_3(j_2\cdots j_aj_1) - \delta_{j_1, j_a}\mu_3(J) \;\; & ( i = j_1),\\
\mu_3(j_2\cdots j_a j_1)\delta_{j_1,j_a} - \mu_3(J) \;\;&  (i = j_a),\\
0 \; \; & \mbox{(otherwise)},
\end{array}\right.\\
& = 0 \;\; (\mbox{by the assumption}).
\end{array}
\leqno{(2.4.5)}
$$
Next, we prove that for any subsequence $J$ of $I'$, $1\leq i \leq r$ and $e = \pm 1$,
$$ \mu_3(J; (x_i^{{\rm N}\frak{p}_i -1})^e)  = 0.  \leqno{(2.4.6)}$$
Suppose $i \in I$. Then  ${\rm N}\frak{p}_i \equiv 1$ mod $m_I$ and so  we write  ${\rm N}\frak{p}_i  - 1= m_I q$. By the assumption, we have
$$ \begin{array}{ll}
\Theta_{r}((x_i^{{\rm N}\frak{p}_i -1})^e ) & = (1 + X_i)^{em_I q}\\
& = \displaystyle{ \left( \sum_{j=0}^{m_I} {m_I \choose j}  X_i^j \right)^{eq} }\\
 & =  1 + ( \mbox{term of degree} \geq |I|),
  \end{array}$$
  from which (2.4.6) follows.\\
Suppose $ i \notin I$. Then $\Theta_{r}((x_i^{{\rm N}\frak{p}_i -1})^e )$ does not contain the term of $X_{J}$, and hence $\mu_3(J; (x_i^{{\rm N}\frak{p}_i -1})^e ) = 0$.
So (2.4.6) is proved.\\
By (2.1.1), (2.4.3) and (2.4.6), we have 
$$ \begin{array}{ll} \mu_3( J; (x_i^{{\rm N}\frak{p}_i -1} [x_i,y_i] )^e)  & =  \displaystyle{\sum_{J', J''} \mu_3(J'; (x_i^{{\rm N}\frak{p}_i -1})^e ) \mu_3(J''; [x_i,y_i]^e)}\\
 & = 0, 
 \end{array}$$
where $J = J'J''$ if $e=1$ and $J = J''J'$ if $e=-1$. Thus (2.4.2) is proved. \\
By the same argument as in the proof of (ii), we have, for any $f \in \frak{F}_{r}$,
$$ \mu_3(J; f(x_i^{{\rm N}\frak{p}_i -1} [x_i,y_i] )^ef^{-1}) = \mu_3(J; (x_i^{{\rm N}\frak{p}_i -1} [x_i,y_i] )^e) = 0.  \leqno{(2.4.7)}$$
Finally, by (2.1.1) and (2.4.7), we have
$$ \begin{array}{ll}
\mu_3(I'; f(x_i^{{\rm N}\frak{p}_i -1} [x_i,y_i] )^ef^{-1}y_{i_n}) & = \displaystyle{\sum_{I' = JK} \mu_3(J; f(x_i^{{\rm N}\frak{p}_i -1} [x_i,y_i] )^ef^{-1}) \mu_3(K;y_{i_n})}\\
& = \mu_3(I).
\end{array}
$$
By the argument similar to the above, we can prove that
$$ \mu_3(I'; y_{i_n}f(x_i^{{\rm N}\frak{p}_i -1} [x_i,y_i] )^ef^{-1})  =  \mu_3(I).$$
Hence the assertion (iii) is proved. $\;\; \Box$
\\
\\
{\bf Remark 2.5.} 
For a multi-index $I$, put $S'_0=\{\frak{p}_i|i \in I\} \subset S_0$, $r'=\#S'_0$ and $S'=S'_0 \cup \{\frak{p}_0\} \subset S$. 
Then the commutative diagram 
\[
\begin{array}{ccccc}
\frak{G}_S & \longrightarrow & \frak{G}_{S'} &  & \\
\uparrow &  & \uparrow &  & \\
\frak{F}_r  & \longrightarrow & \frak{F}_{r'} & : & x_j \mapsto 1 \ (j \not\in I) \\
\downarrow &  & \downarrow &  & \\
\mathbb F_3\langle\langle X_1,\ldots,X_r \rangle\rangle  & \longrightarrow & \mathbb F_3\langle\langle X_i\ (i \in I) \rangle\rangle & : & X_j \mapsto 0 \ (j \not\in I) \\
\end{array}
\]
yields that $\mu_3^{S}(I)=\mu_3^{S'}(I)$. This reduces our consideration on mod $3$ Milnor invariants to the case where $S_0 =\{ \frak{p}_{i_1},\dots , \frak{p}_{i_n} \}$ and $I = (i_1\cdots i_n)$.
\\

We end this section by giving the following \\
\\
{\bf Definition 2.6.} Let $\frak{p}_{1},\dots , \frak{p}_{n-1}$ be distinct $n-1$ primes in $S_k^{1\, \mbox{{\scriptsize mod}}\, 9}$.  We call an extension $K$ of $k$ a {\em R\'{e}dei type $H_n(\mathbb{F}_3)$-extension for $\{ \frak{p}_{1},\dots , \frak{p}_{n-1}\}$} if $K$ is a Galois extension of $k$ in which only primes $\frak{p}_{1},\dots , \frak{p}_{n-1}$ are ramified with ramification index of each $\frak{p}_{j}$ being $3$ and whose Galois group is isomorphic to the group $H_n(\mathbb{F}_3)$ consisting of $n$ by $n$ upper-triangular unipotent matrices over $\mathbb{F}_3$:
$$
H_n(\mathbb{F}_3) := \left\{ \left(
\begin{array}{ccccc}
	 1 & * & * & \cdots & * \\
	 0 & 1 & * & \cdots & * \\
	 \vdots & \ddots & \ddots & \ddots & \vdots \\
	 \vdots & & \ddots & 1 &  * \\
	 0 & \cdots & \cdots & 0 & 1 \\
\end{array}
\right) \; \vert \; * \in \mathbb{F}_3 \right\}.
$$
\vspace{.5cm}
\\
As we shall see in Sections 3 and 4, the independence of $\mu_3(1\cdots n)$ on the choice of $\frak{p}_0$ will be deduced from the existence of the R\'{e}dei type $H_n(\mathbb{F}_3)$-extension for $\{ \frak{p}_{1},\dots , \frak{p}_{n-1}\}$, when $n =2$ or $3$.\\

\begin{center}
{\bf 3. Mod $3$ linking numbers and cubic residue symbols}
\end{center}

In this section, we prove that the mod $3$ Milnor number  $\mu_3(12)$ is an invariant (mod $3$ linking number) determined by  primes $\frak{p}_1, \frak{p}_2 \in  S_k^{1\, \mbox{{\scriptsize mod}}\, 9 }$, by showing that there exists (uniquely) a R\'{e}dei type $H_2(\mathbb{F}_3)$-extension $\frak{K}_{\{\frak{p}_1\}}$ of $k$ for $\{\frak{p}_1\}$. We then show its relation with the cubic residue symbol by constructing $\frak{K}_{\{\frak{p}_1\}}$ concretely.  \\

We start to show the following
\\
\\
{\bf Theorem 3.1.} {\em For $\frak{p} \in S_k^{1\, \mbox{{\scriptsize mod}}\, 9 }$, there exists uniquely a Redei type $H_2(\mathbb{F}_3)$-extension $\frak{K}_{\{\frak{p}\}}$ of $k$ for $\{\frak{p}\}$. }\\
\\
{\em Proof.}  Note $H_2(\mathbb{F}_3) = \mathbb{Z}/3\mathbb{Z}$. By Proposition 1.8, we have $d(G_{\{ \frak{p} \}}) = 1$, which implies the assertion. $\;\; \Box$\\
\\
Let $S_0 = \{ \frak{p}_1, \frak{p}_2 \}$ for the simplicity of the notation in Section 2. By using Theorem 3.1, we see the following \\
\\
{\bf Theorem 3.2.} {\em The mod $3$ Milnor number $\mu_3(12)$ is independent of the choice of $\frak{p}_0$ and so an invariant of the ordered pair $(\frak{p}_1, \frak{p}_2)$.}\\
\\
{\em Proof.} Take a prime $\frak{p}_0 \in  S_k^{4,7 \, \mbox{{\scriptsize mod}}\, 9}$ and let  $S := S_0 \cup \{ \frak{p}_0\} = \{ \frak{p}_0, \frak{p}_1, \frak{p}_1 \}$. 
We use the same notations as in Sections 1 and 2. By the definition of $\mu_3(12)$,  we have 
$$ \sigma_2 \equiv \tau_1^{\mu_3(12)} \; \mbox{mod} \; \frak{G}_S^{(2)}. \leqno{(3.2.1)}$$ 
There exists uniquely an $H_2(\mathbb{F}_3)$-extension $\frak{K}_{\{\frak{p}_1\}}$ of $k$ for $\{\frak{p}_1\}$ by Theorem 3.1 and the Galois group ${\rm Gal}(\frak{K}_{\{\frak{p}_1\}}/k) \simeq \mathbb{Z}/3\mathbb{Z}$ is generated by $\tau_1|_{\frak{K}_{\{\frak{p}_1\}}}$. Since $\frak{K}_{\{\frak{p}_1\}} \subset k_S$,  we have, by (3.2.1), 
$$ \sigma_2|_{\frak{K}_{\{\frak{p}_1\}}} = (\tau_1|_{\frak{K}_{\{\frak{p}_1\}}})^{\mu_3(12)},$$
which implies that $\mu_3(12)$ is independent of the choice of $\frak{p}_0$ and so an invariant of $S_0$. $\;\; \Box$\\
\\
We call $\mu_3(12)$ the {\it mod $3$ linking number} of $\frak{p}_1$ and $\frak{p}_2$ and denote it by ${\rm lk}_3(\frak{p}_1,\frak{p}_2)$.\\

Next, we show the relation between the mod $3$ linking number and the cubic residue symbol by constructing concretely $\frak{K}_{\{ \frak{p}\}}$ in Theorem 3.1. For this, we  recall a well-known fact on the ramification in a Kummer extension. \\
\\
{\bf Lemma 3.4} ([B; Lemma 5, Lemma 6]). {\em Let $F$ be a finite algebraic number field containing a primitive $l$-th root of unity, where $l$ is a prime number. Let  $L$ be a cyclic extension of degree $l$ over $F$ so that $L = F(\sqrt[l]{a})$ for some $a \in F^{\times}$. Let $\frak{p}$ be a finite prime of $F$ which is not lying over $l$, and write $(a) = \frak{p}^{v_{\frak{p}}(a)}\frak{q}$ with $\frak{q}$ being an ideal prime to $\frak{p}$.  Then we have the follwings.}\\
(1) {\em If $v_{\frak{p}}(a)$ is prime to $l$,  then $\frak{p}$ is totally ramified in $L/F$.}\\
(2) {\em If $v_{\frak{p}}(a)$ is divisible by $l$, then $\frak{p}$ is unramified in $L/F$.}\\
\\
{\bf Theorem 3.5.} {\em Let $\frak{p} \in  S_k^{1\, \mbox{{\scriptsize mod}}\, 9}$ and let $\pi$ the prime element in $\frak{p}$ satisfying $\pi \equiv 1 \; \mbox{mod} \; (3\sqrt{-3})$ as in Lemma 1.1. Then $k(\sqrt[3]{\pi})$ is a cyclic extension of degree $3$ over $k$ in which only $\frak{p}$ is ramified. Hence we have}
$$ \frak{K}_{\{ \frak{p} \}} = k(\sqrt[3]{\pi}).$$
\\
{\em Proof.} It is easy to see that $k(\sqrt[3]{\pi})$ is a cyclic Kummer extension of degree $3$ over $k$. So it suffices to show the assertion on the ramification. Since $(\pi) = \frak{p}$, $\frak{p}$ is totally ramified by Lemma 3.4 (1). Let $\lambda := \frac{\sqrt[3]{\pi}-1}{\sqrt{-3}}$. Since $\lambda$ satisfies 
$\lambda^3 - \sqrt{-3} \lambda^2 - \lambda - \frac{1 - \pi}{3\sqrt{-3}} = 0$ and $\frac{1 - \pi}{3\sqrt{-3}} \in {\cal O}_k$ by Lemma 1.1, we find $\lambda \in {\cal O}_K$. 
 The relative discriminant of $\lambda$ in $k(\sqrt[3]{\pi})/k$ is computed as 
$$
d(\lambda, K/k) = 
\left|
\begin{array}{ccc}
1 & \lambda^{(1)} & (\lambda^{(1)})^2 \\
1 & \lambda^{(2)} & (\lambda^{(2)})^2 \\
1 & \lambda^{(3)} & (\lambda^{(3)})^2 \\
\end{array}
\right|^2
= -\frac{\pi^2}{27}
\left|
\begin{array}{ccc}
1 & 1 & 1 \\
1 & \zeta_3 & \zeta_3^2 \\
1 & \zeta_3^2 & \zeta_3 \\
\end{array}
\right|^2
= \pi^2,
$$
where $\lambda^{(1)} := \lambda$,  $\lambda^{(2)} := (\zeta_3 \sqrt[3]{\pi} - 1)/\sqrt{-3}$ and $\lambda^{(2)} := (\zeta_3^2 \sqrt[3]{\pi} -1)/\sqrt{-3}$. Hence $k(\sqrt[3]{\pi})/k$ is unramified outside $\frak{p}$. The second assertion follows from Theorem 3.1. $\;\; \Box$\\
\\
Recall that for a finite prime $\frak{p} = (\pi)$ of $k$ and $a \in k_{\frak{p}}^{\times}$ such that $v_{\frak{p}}(a) \equiv 0$ mod $3$, the {\em cubic residue symbol} is defined by 
$$ \left(  \frac{a}{\pi} \right)_3 = \frac{\sigma(\sqrt[3]{a})}{\sqrt[3]{a}},$$
where $\sigma$ is the Frobenius automorphism of the unramified extension $k_{\frak{p}}(\sqrt[3]{a})/k_{\frak{p}}$. We note that
$$ \left(  \frac{a}{\pi} \right)_3 = 1 \; \Longleftrightarrow  \; a \in (k_{\frak{p}}^{\times})^3.$$
\\
{\bf Theorem 3.6.} {\em Let $\frak{p}_1$ and $\frak{p}_2$ be distinct primes in $S_k^{1 \, \mbox{{\scriptsize mod}}\, 9}$ and let $\pi_i$ be the prime element in $\frak{p}_i$ satisfying $\pi_i \equiv 1 \; \mbox{mod} \; (3\sqrt{-3})$ as in Lemma 1.1 for $i=1,2$. Then we have }
$$   \zeta_3^{{\rm lk}_3(\frak{p}_1,\frak{p}_2)} = \left(  \frac{\pi_1}{\pi_2} \right)_3.$$
\\
{\em Proof.} This is just a refinement of the proof of Theorem 3.2. Take a prime $\frak{p}_0 \in  S_k^{4,7 \, \mbox{{\scriptsize mod}}\, 9}$ and let  $S := S_0 \cup \{ \frak{p}_0\} = \{ \frak{p}_0, \frak{p}_1, \frak{p}_2  \}$. By the definition of $\mu_3(12) =: {\rm lk}_3(\frak{p}_1,\frak{p}_2)$, we have
$$ \sigma_2 \equiv \tau_1^{{\rm lk}_3(\frak{p}_1,\frak{p}_2)} \; \mbox{mod} \; \frak{G}_S^{(2)}. \leqno{(3.6.1)}$$ \\
Noting that $\frak{K}_{\{ \frak{p}_1 \}} = k(\sqrt[3]{\pi_1}) \subset k_S$, we have, by (1.2), 
$$ \tau_1(\sqrt[3]{\pi_1}) = \zeta_3 \sqrt[3]{\pi_1}. \leqno{(3.6.2)}$$
By the definition of the cubic residue symbol as in the above, we have
$$ \left(  \frac{\pi_1}{\pi_2}  \right)_3 = \frac{\sigma_2(\sqrt[3]{\pi_1})}{\sqrt[3]{\pi_1}}. \leqno{(3.6.3)}$$
By (3.6.1), (3.6.2) and (3.6.3), we obtain the assertion. $\;\; \Box$\\

\begin{center}
{\bf 4. Mod $3$ triple Milnor invariants}
\end{center}

In this section,  we give a condition for a R\'{e}dei type $H_3(\mathbb{F}_3)$-extension $\frak{K}_{\{ \frak{p}_1, \frak{p}_2\}}$ of $k$ for $\{\frak{p}_1, \frak{p}_2\}$ to exist, and show that the condition is satisfied if $\frak{p}_1$ and $\frak{p}_2$ are generated by prime numbers. We then prove that the mod $3$ triple Milnor number $\mu_3(123)$ is an invariant of the ordered triple $(\frak{p}_1, \frak{p}_2, \frak{p}_3)$ of primes  in $S_k^{1 \, \mbox{{\scriptsize mod}}\, 9}$, under $\mu_3(ij) = 0$ ($1\leq i \neq j \leq 3$), if there exists (uniquely) a R\'{e}dei type $H_3(\mathbb{F}_3)$-extension $\frak{K}_{\{ \frak{p}_1, \frak{p}_2\}}$ of $k$ for $\{\frak{p}_1, \frak{p}_2\}$.
\\

Let $\frak{p}_1$ and $\frak{p}_2$ be distinct primes in $S_k^{1 \, \mbox{{\scriptsize mod}}\, 9}$ and let  $\pi_i$ be the unique prime element in $\frak{p}_i$ such that $\pi_i \equiv \; 1 \; \mbox{mod}\; (3\sqrt{-3})$ as in Lemma 1.1. 
As in Section 3, mod $3$ linking numbers $\mu_3(ij)$ are well defined for $1\leq i \neq j \leq 3$. We set 
$$ K_i := k(\sqrt[3]{\pi_i}) = \frak{K}_{\{ \frak{p}_i \}} \; (i=1,2) \; \mbox{and}\;  K_{12}=k(\sqrt[3]{\pi_1 \pi_2}). $$ 
First, we prove the following.\\
\\
{\bf Theorem 4.1.} {\em 
There exists a R\'{e}dei type  $H_3(\mathbb{F}_3)$-extension $\frak{K}_{\{\frak{p}_1, \frak{p}_2\}}$ of $k$ for $\{\frak{p}_1, \frak{p}_2\}$ 
if and only if the class number $h_{K_{12}}$ of $K_{12}$ is divisible by $9$. Moreover, $\frak{K}_{\{\frak{p}_1, \frak{p}_2\}}$ is unique if it exists. }
\\
\\
{\em Proof.} 
Let $L$ be the Hilbert $3$-class field of $K_{12}$. 
The genus $3$-class field of $K_{12}/k$ is the maximal subfield $L^g$ of $L$ which is abelian over $k$. 
Then $K_1K_2 \subset L^g$. 
Since the class number $h_k$ is one, 
$\mathrm{Gal}(L^g/k)$ is generated by the inertia subgroups $I_i$ of $\mathfrak{p}_i$. 
Then $[L^g:k] \le |I_1| |I_2| =9$ and hence $L^g=K_1K_2$. 
In particular, the maximal abelian quotient of $G=\mathrm{Gal}(L/k)$ is 
$G/[G,G] \simeq \mathrm{Gal}(K_1K_2/k) \simeq (\mathbb{Z}/3\mathbb{Z})^2$. 
Let $\mathfrak{K}$ be the fixed field of $[G,G]^3[[G,G],G]$. 

If part: Suppose that $h_{K_{12}}$ is divisible by $9$, i.e., $[L:K_{12}] \ge 9$. Then $G$ is nonabelian. 
Hence 
$\mathrm{Gal}(\mathfrak{K}/k) \simeq G/[G,G]^3[[G,G],G]$ 
is a nonabelian $3$-group of order $27$. 
Since $\mathfrak{K}/K_1K_2$ is unramified, 
neither of $\mathfrak{K}/K_1$ and $\mathfrak{K}/K_2$ is cyclic extensions, 
i.e., 
the noncyclic maximal subgroup of $\mathrm{Gal}(\mathfrak{K}/k)$ 
is not unique. 
This implies that $\mathrm{Gal}(\mathfrak{K}/k) \simeq H_3(\mathbb F_3)$. 
Therefore $\mathfrak{K}$ is a R\'{e}dei type  $H_3(\mathbb{F}_3)$-extension of $k$ for  $\{ \frak{p}_1, \frak{p}_2\}$. 

Only-if part and uniqueness: 
Suppose the existence of $\frak{K}_{\{\frak{p}_1, \frak{p}_2\}}$. 
Since $d(\frak{G}_{\{\frak{p}_1, \frak{p}_2\}})=2$ by Proposition 1.8, 
$K_1K_2$ is the unique $(\mathbb Z/3\mathbb Z)^2$-extension of $k$ unramified outside $\{\frak{p}_1, \frak{p}_2\}$. 
Since $H_3(\mathbb{F}_3)/[H_3(\mathbb{F}_3),H_3(\mathbb{F}_3)] \simeq (\mathbb Z/3\mathbb Z)^2$, we have $K_1K_2 \subset \frak{K}_{\{\frak{p}_1, \frak{p}_2\}}$. 
Since $\mathfrak{p}_1$ and $\mathfrak{p}_2$ ramify in $K_{12}/k$, 
$\frak{K}_{\{\frak{p}_1, \frak{p}_2\}}/K_{12}$ is an unramified Galois extension of degree $9$. 
Hence $\frak{K}_{\{\frak{p}_1, \frak{p}_2\}} \subset L$. By class field theory, $h_{K_{12}}$ is divisible by $9$. 
Since 
\[
[G,G]^3[[G,G],G] \subset \mathrm{Ker}(G \stackrel{\mathrm{restr.}}{\longrightarrow} \mathrm{Gal}(\frak{K}_{\{\frak{p}_1, \frak{p}_2\}}/k)) , 
\]
$\frak{K}_{\{\frak{p}_1, \frak{p}_2\}}$ coincides with the fixed field $\mathfrak{K}$ of $[G,G]^3[[G,G],G]$, which is uniquely defined. 
$\;\; \Box$\\
\\
{\bf Proposition 4.2.} {\em 
If $\frak{p}_1=(p_1)$ and $\frak{p}_2=(p_2)$ with prime numbers $p_1$ and $p_2$, then the class number of $K_{12}$ is divisible by $9$ and hence  there exists uniquely a R\'{e}dei type  $H_3(\mathbb{F}_3)$-extension $\frak{K}_{\{\frak{p}_1, \frak{p}_2\}}$ of $k$ for $\{\frak{p}_1, \frak{p}_2\}$.}\\
\\
{\em Proof.} 
Since $p_i$ is inert in $k/\mathbb Q$ and $p_i^2 \equiv 1 \pmod{9}$, 
we have $p_i \equiv -1 \pmod{9}$. 
Then $\pi_i=-p_i$. Hence $K_i=k(\sqrt[3]{p_i})$ and $K_{12}=k(\sqrt[3]{p_1p_2})$. 
Note that $K_1K_2/K_{12}$ is an unramified cubic cyclic extension, 
and that 
\[
\mathrm{Gal}(K_1K_2/\mathbb Q(\sqrt[3]{p_1p_2})) \simeq \mathrm{Gal}(K_i/\mathbb Q) \simeq S_3 . 
\]
On the other hand, there is an unramified cubic cyclic extension $F/\mathbb Q(\sqrt[3]{p_1p_2})$ by [H]. 
Since 
\[
\mathrm{Gal}(K_{12}F/\mathbb Q(\sqrt[3]{p_1p_2})) \simeq \mathbb Z/6\mathbb Z, 
\]
$K_{12}F/K_{12}$ is another unramified cubic cyclic extension. 
Therefore $K_1K_2F/K_{12}$ is an unramified $(\mathbb Z/3\mathbb Z)^2$-extension. Hence  $h_{K_{12}}$ is divisible by $9$ by class field theory.  By Theorem 4.1, the latter assertion follows.
$\;\; \Box$\\

Let $S_0 = \{ \frak{p}_1, \frak{p}_2, \frak{p}_3 \}$ for the simplicity of notation in Section 2. Let $\pi_i$ be the unique prime element  in $\frak{p}_i$ satisfying $ \pi_i \equiv 1 \; \mbox{mod} \;  (3\sqrt{-3}) \;  (1\leq i \leq 3)$ as in Lemma 1.1. We assume that
$$ \mu_3(ij) = {\rm lk}_3(\frak{p}_i,\frak{p}_j) = 0 \;\;\;\; (1\leq i\neq j \leq 3),   \leqno{(4.3)}$$
which is equivalent, by Theorem 3.6, to the condition
$$ \left( \frac{\pi_i}{\pi_j} \right)_3 = 1 \;\;\;\; (1\leq i\neq j \leq 3). $$
Now, we prove the following \\
\\
{\bf Theorem  4.4.} {\em 
If there exists a R\'{e}dei type $H_3(\mathbb{F}_3)$-extension $\frak{K}_{\{\frak{p}_1, \frak{p}_2\}}$ of $k$ for $\{\frak{p}_1, \frak{p}_2\}$, then the
mod $3$ Milnor number $\mu_3(123)$ is independent of the choice of $\frak{p}_0$ and so an invariant of the ordered triple $(\frak{p}_1, \frak{p}_2, \frak{p}_3)$.}\\
\\
{\em Proof.} Choose a prime $\frak{p}_0  \in S_k^{4, 7\, \mbox{{\scriptsize mod}}\, 9}$ and set $S := \{ \frak{p}_0, \frak{p}_1, \frak{p}_2, \frak{p}_3\}$. By Theorem 1.10, the Galois group $\frak{G}_S$ of the maximal pro-$3$ extension $k_S$ over $k$ unramified outside $S$ has the minimal presentation:
$$ \begin{array}{ll} \frak{G}_{S} & = \langle x_1, x_2, x_3  \mid  x_{1}^{{\rm N}\frak{p}_1- 1}[x_1 ,y_1] =   x_{2}^{{\rm N}\frak{p}_2- 1}[x_2 ,y_2] = x_{3}^{{\rm N}\frak{p}_3- 1}[x_3 ,y_3] = 1 \rangle \\
& = \frak{F}_3/\frak{N}_S.
\end{array}$$
Here $x_i$ is the letter representing a monodromy $\tau_i$ over $\frak{p}_i$ in $k_S/k$ and $\frak{F}_3$ is the free pro-$3$ group generated by $x_1, x_2$ and $x_3$. 
The free pro-$3$ word $y_i \in \frak{F}_3$  represents a Frobenius automorphism $\sigma_i$ over $\frak{p}_i$ in $k_S/k$ and $\frak{N}_S$ denotes the closed subgroup of
$\frak{F}_3$ generated normally by $x_{i}^{{\rm N}\frak{p}_i - 1}[x_i ,y_i]$ for $1\leq i \leq 3$. By (2.1.3) and the assumption (4.3), we have, for $1\leq i \leq 3$, 
$$ y_i \in \frak{F}_r^{(2)} = \frak{F}_3^3[\frak{F}_3,\frak{F}_3].  \leqno{(4.4.1)}$$
Let $\frak{F}(s_i, t_i)$ be the free pro-$3$ group on words $s_i,t_i$ $(1\leq i \leq 3)$ and let $\tilde{k}$ denote the maximal pro-$3$ extension of $k$. By [Ko; Satz 6.11], we have the following commutative diagram
$$ 
\setlength{\unitlength}{0.7mm}
\begin{picture}(100,50)
\put(-10,0){\makebox(20,10)[r]{$\frak{F}_3$}}
\put(-5,30){\makebox(20,10)[r]{$\frak{F}(s_i,t_i)$}}
\put(39,0){\makebox(20,10){$\frak{G}_S$}}
\put(39,30){\makebox(20,10){$\frak{G}_{\frak{p}_i}$}}
\put(80,0){\makebox(20,10)[l]{$\mathrm{Gal}(\mathfrak{K}_{\{ \frak{p}_1,\frak{p}_2\}}/k)$}}
\put(80,30){\makebox(20,10)[l]{${\rm Gal}(\tilde{k}/k)$.}}
\put(22,5){\vector(1,0){18}}
\put(60,5){\vector(1,0){18}}
\put(22,35){\vector(1,0){18}}
\put(60,35){\vector(1,0){18}}
\put(5,27){\vector(0,-1){15}}
\put(47,27){\vector(0,-1){15}}
\put(88,27){\vector(0,-1){15}}
\put(80,28){\vector(-3,-2){25}}
\put(65,0){$\Phi_S$}
\put(90,17){$\Phi$}
\put(65,38){$\varphi_{\frak{p}_i}$}
\put(32,18){$\varphi_{\frak{p}_i, S}$}
\put(25,38){${\scriptstyle t_i \mapsto \tau_i}$}
\put(25,42){${\scriptstyle s_i \mapsto \sigma_i}$}
\put(-10,15){${\scriptstyle t_i \mapsto x_i}$}
\put(-10,20){${\scriptstyle s_i \mapsto y_i}$}
\end{picture}
$$
Here $\varphi_{\frak{p}_i, S}$ is the map in (1.4) and $\varphi_{\frak{p}_i}$ is also the map induced by the embedding $\tilde{k} \hookrightarrow \tilde{k}_{\frak{p}_i}$. 
The maps $\Phi_S$ and $\Phi$ are natural quotient homomorphisms. We set, for $1\leq i \leq 3$, 
$$
\gamma_i :=\Phi(\varphi_{\mathfrak{p}_i}(\tau_i))=\Phi_S(x_i\mathfrak{N}_S), \;
\eta_i :=\Phi(\varphi_{\mathfrak{p}_i}(\sigma_i))=\Phi_S(y_i\mathfrak{N}_S) \;
\in \mathrm{Gal}(\mathfrak{K}_{\{ \frak{p}_1, \frak{p}_2 \}}/k).
$$
We note that $\gamma_i$ and $\eta_i$ are independent of the choice of $\mathfrak{p}_0$.  Let $\frak{V}$ denote the closed subgroup of $\mathfrak{F}_3$ generated normally by 
 $x_1^3,x_2^3$ and $x_3$. Then we have $x_{3}^{{\rm N}\frak{p}_3 - 1}[x_3 ,y_3] \in \frak{V}$ obviously and $x_{i}^{{\rm N}\frak{p}_i - 1}[x_i ,y_i] \in \frak{V}[[\frak{F}_3,\frak{F}_3],\frak{F}_3]$ for $i=1,2$ by (4.4.1). Therefore we have $\mathfrak{N}_S \subset \frak{V}[[\mathfrak{F}_3,\mathfrak{F}_3],\mathfrak{F}_3]$.
We set 
$$\mathfrak{M} :=  \frak{V}\frak{F}_3^{(3)} =  \frak{V}\mathfrak{F}_3^3[[\mathfrak{F}_3,\mathfrak{F}_3],\mathfrak{F}_3],$$
where $\frak{F}_3^{(3)}$ is the 3rd term of the Zassenhaus filtration of $\frak{F}_3$ (cf. (2.1.3)). Since $H_3(\mathbb F_3)$ is two-step nilpotent and has the exponent $3$, 
the surjective homorophsim $\Phi_S$ factors through 
$\mathfrak{F}_3/\mathfrak{M}$. 
Since $\mathfrak{F}_3/[\mathfrak{F}_3,\mathfrak{F}_3]\mathfrak{M}$
is generated by $x_1, x_2 \bmod [\mathfrak{F}_r,\mathfrak{F}_r]\mathfrak{M}$ and the commutator subgroup of $\mathfrak{F}_3/\mathfrak{M}$ 
is a cyclic group generated by $[x_1,x_2] \bmod \mathfrak{M}$, 
we see that 
$$
\#(\mathfrak{F}_3/\mathfrak{M}) \leq  27 . 
$$
Therefore $\Phi_S$ induces the isomorphism
$$
\mathfrak{F}_3/\mathfrak{M} \simeq \mathrm{Gal}(\mathfrak{K}_{\{ \frak{p}_1,\frak{p}_2\}}/k),  \leqno{(4.4.2)}
$$
where $x_i \; \mbox{mod}\; \frak{M}$ is sent to $\gamma_i$ for $i=1,2$ and to ${\rm id}$ for $i=3$, and we note that the commutator subgroup of $\mathrm{Gal}(\mathfrak{K}_{\{ \frak{p}_1,\frak{p}_2\}}/k)$ 
is a cyclic group generated by $[\gamma_1,\gamma_2]$.  Since $y_3 \in \frak{V}\mathfrak{F}_3^3[\mathfrak{F}_3,\mathfrak{F}_3]$ by (4.4.1), 
there exists uniquely $m(123) \in \mathbb F_3$ such that 
$$
y_3 \equiv [x_1,x_2]^{m(123)} \mod{\mathfrak{M}}. \leqno{(4.4.3)}
$$
By (2.1.1) and (2.1.3), we have $\mu_3(I; \frak{M})=0$ for $I \in \{(1),(2),(12)\}$ and 
$$
\mu_3(123) := \mu_3((12); y_3) = m(123) \in \mathbb F_3. \leqno{(4.4.4)}
$$
By (4.4.2), (4.4.3) and (4.4.4), we have 
$$
\eta_3 = [\gamma_1,\gamma_2]^{\mu_3(123)} \in \mathrm{Gal}(\mathfrak{K}_{\{ \frak{p}_1, \frak{p}_2 \}}/k), 
$$
which implies that $\mu_3(123)$ is independent of the choice of $\frak{p}_0$ and so an invariant of $S_0$.  $\;\; \Box$\\
\\
Composing (4.4.2) with the isomorphism  ${\rm Gal}(\frak{K}_{\{ \frak{p}_1, \frak{p}_2\}}/k) \simeq H_3(\mathbb{F}_3) $ given by the  correspondence 
$$ \gamma_1 \mapsto \left( \begin{array}{ccc} 1 & 1 & 0 \\ 0 & 1 & 0 \\ 0 & 0 & 1 \end{array}\right), \; \; \gamma_2 \mapsto \left( \begin{array}{ccc} 1 & 0 & 0 \\ 0 & 1 & 1\\ 0 & 0 & 1 \end{array}\right),$$
we have the isomorphism
$$ \rho : \frak{F}_3/\frak{M} \stackrel{\sim}{\rightarrow} H_3(\mathbb{F}_3).$$
We note that $\rho$ is given by
$$ \rho(f \; \mbox{mod} \; \frak{M}) = \left( \begin{array}{ccc} 1 & \mu_3((1);f) & \mu_3((12);f) \\ 0 & 1 & \mu_3((2);f) \\ 0 & 0 & 1 \end{array}\right)$$
and we have
$$ \rho(y_3 \; \mbox{mod} \; \frak{M}) = \left( \begin{array}{ccc} 1 & 0 & \mu_3(123) \\ 0 & 1 & 0 \\ 0 & 0 & 1 \end{array}\right).$$
We remark that a map analogous to $\rho$ was considered for the study of Milnor invariants in link theory ([Mu]).\\
\\
{\bf Corollary 4.5.} {\em If $\frak{p}_1=(p_1)$ and $\frak{p}_2=(p_2)$ with prime numbers $p_1$ and $p_2$, then the
mod $3$ Milnor number $\mu_3(123)$ is an invariant of the ordered triple $(\frak{p}_1, \frak{p}_2, \frak{p}_3)$.}\\
\\
{\em Proof.} This follows from Proposition 4.2 and Theorem 4.4. $\;\; \Box$\\
\\
There is a case where no R\'{e}dei type $H_3(\mathbb{F}_3)$-extension exists, and $\mu_3^S(123)$ depends on the choice of $\frak{p}_0$ as follows. 
\\
\\
{\bf Example 4.6.} The following examples are calculated by [Pari-gp]. 
Put $\pi_1=-17$, $\pi_2=1+9\zeta_3$ and $\pi_3=7+12\zeta_3$. 
Then $\pi_i \equiv 1 \pmod{3\sqrt{-3}}$ for $1 \le i \le 3$, and ${\rm N}\frak{p}_2=73$, ${\rm N}\frak{p}_3=109$. 
Moreover, the condition (4.3) is satisfied. 
Then $\frak{p}_3$ splits completely in $K_1K_2/k$. 
Put 
\[
K'_{12}=k(\sqrt[3]{\pi_1^2\pi_2}) \simeq \mathbb Q[x]/(x^6+2023x^3+6097033) . 
\]
Since $h_{K'_{12}} = 3$, 
$K_1K_2$ is the maximal unramified $3$-extension of $K'_{12}$, 
in particular $K_1K_2$ has no unramified $3$-extension. 
Hence $K_1K_2$ is the Hilbert $3$-class field of $K_{12}$, 
and $h_{K_{12}}$ is not divisible by $9$. 
On the other hand, $\overline{\frak{p}}_3=(\overline{\pi}_3)=(7+12\zeta_3^{-1})$ is inert in $K_2/k$ and splits in $K_1/k$. 
Then $\overline{\frak{p}}_3$ is also inert in $K'_{12}/k$, 
and hence a prime $\frak{P}_3$ of $K'_{12}$ lying over $\frak{p}_3$ is characterized as a degree one prime of $K'_{12}$ lying over $109$. 
Note that $\frak{P}_3$ splits in $K_1K_2/K'_{12}$. 

Put $\frak{p}_0=(1+3\zeta_3) \ni 7$ or $\frak{p}_0=(-6-7\zeta_3) \ni 43$. 
Then ${\rm N}\frak{p}_0 \equiv 7 \pmod{9}$, and $\frak{p}_0$ is inert in $K'_{12}/k$. 
Moreover, the ray class group $Cl_{\frak{p}_0}(K'_{12})$ of $K'_{12}$ modulo $\frak{p}_0{\cal O}_{K'_{12}}$ has the $3$-rank $2$. 
Let $\frak{K}/K'_{12}$ be the maximal elementary abelian $3$-extension unramified outside $\frak{p}_0$. 
Then $\mathrm{Gal}(\frak{K}/K'_{12}) \simeq Cl_{\frak{p}_0}(K'_{12})/3 \simeq (\mathbb Z/3\mathbb Z)^2$, 
and $\frak{K}/k$ is also a Galois extension. 
Since $\frak{K}/K_i$ is not cyclic, i.e., the noncyclic maximal subgroup of $\mathrm{Gal}(\frak{K}/k)$ is not unique, $\frak{K}/k$ is a $H_3(\mathbb F_3)$-extension unramified outside $\{\frak{p}_0,\frak{p}_1,\frak{p}_2\}$. 
Put $S=\{\frak{p}_0,\frak{p}_1,\frak{p}_2,\frak{p}_3\}$. 
Then the map $\Phi_S : \frak{G}_S \rightarrow \mathrm{Gal}(\frak{K}/k)$ induces an isomorphism $\frak{F}_3/\frak{M} \simeq \mathrm{Gal}(\frak{K}/k)$ similar to (4.4.2), 
and we have 
\[
y_3 \equiv [x_1,x_2]^{\mu_3^S(123)} \mod{\frak{M}} 
\]
as in the proof of Theorem 4.4. Let $[\frak{P}_3]$ be the class of $\frak{P}_3$ in $Cl_{\frak{p}_0}(K'_{12})$. 
\\
$\bullet$ If $7 \in \frak{p}_0=(1+3\zeta_3)$, then $[\frak{P}_3] \not\in Cl_{\frak{p}_0}(K'_{12})^3$, i.e., a prime of $K_1K_2$ lying over $\frak{P}_3$ is inert in $\frak{K}/K_1K_2$. This implies that $\mu_3^S(123) \neq 0 \in \mathbb F_3$. 
\\
$\bullet$ If $43 \in \frak{p}_0=(-6-7\zeta_3)$, then $[\frak{P}_3] \in Cl_{\frak{p}_0}(K'_{12})^3$, i.e., $\frak{p}_3$ splits completely in $\frak{K}/k$. This implies that $\mu_3^S(123) = 0 \in \mathbb F_3$. \\

\begin{center}
{\bf 5. Construction of the R\'{e}dei type $H_3(\mathbb{F}_3)$-extension for $\{ \frak{p}_1, \frak{p}_2 \}$}
\end{center}

In this section, we construct concretely the R\'{e}dei type $H_3(\mathbb{F}_3)$-extension $\frak{K}_{\{\frak{p}_1,\frak{p}_2\}}$ of $k$ for $\{ \frak{p}_1, \frak{p}_2\}$,  under a certain assumption (A), in the form analogous to R\'{e}dei's dihedral extension over $\mathbb{Q}$.  We show that the assumption (A)  holds if $\frak{p}_1$ and $\frak{p}_2$ are generated by prime numbers.\\

We start to recall a classical theorem by Iwasawa and the ambiguous class number formula.\\
\\
{\bf Lemma 5.1} ([Iw1; II]). {\em Let $l$ be a prime number and let $F$ be a finite algebraic number field whose class number $h_F$ is prime to $l$. Let $L$ be a Galois extension over $F$ of degree a power of $l$. We assume that there is at most one prime in $F$ which is ramified in $L/F$. Then the class number $h_L$ is also prime to $l$.}\\
\\
{\bf Lemma 5.2} ([L; Theorem 1 (2)], [Y; Lemma 4]). {\em Let $l$ be a prime number and let $L/F$ be a cyclic extension of finite algebraic number fields of degree $l$. Let $A(L/F)$ be the group consisting of classes of ideals $\frak{a}$ of $L$ satisfying $\tau(\frak{a}) = \frak{a}$, where $\tau$ is a generator of the Galois group ${\rm Gal}(L/F)$. Then we have}
$$ \#A(L/F) = \frac{h_F \cdot l^t}{[{\cal O}_F^{\times} : {\rm N}_{L/F}({\cal O}_L^{\times})]},$$
{where $t$ is the number of primes of $F$ which are ramified in $L/F$.}\\
\\
Using Lemma 5.1 and Lemma 5.2, we have the following Propositions.\\
\\
{\bf Proposition 5.3.} {\em Let $\frak{p}$ be a prime in $S_k^{1\, \mbox{{\scriptsize mod}}\, 9}$ and let $\pi$ be the unique prime element in $\frak{p}$ satisfying $\pi \equiv 1 \; \mbox{mod}\; (3\sqrt{-3})$ as in Lemma 1.1. Then the class number of $k(\sqrt[3]{\pi})$ is prime to $3$.}\\
\\
{\em Proof.} This follows from Theorem 3.5 and Lemma 5.1. $\;\; \Box$\\
\\
{\bf Proposition 5.4.} {\em Let $\frak{p}$ and $\pi$ be as in Proposition 5.3 and set $K := k(\sqrt[3]{\pi})$. Then the norm map restricted to the unit groups ${\rm N}_{K/k} : {\cal O}_K^{\times} \rightarrow {\cal O}_k^{\times}$ is surjective.}
\\
\\
{\em Proof.} This follows from Theorem 3.5, Lemma 5.2 and $h_k = 1$. $\;\; \Box$\\
\\
Let $\frak{p}_1$ and $\frak{p}_2$ be distinct finite primes in $S_k^{1\, \mbox{{\scriptsize mod}}\, 9}$. By Lemma 1.1, we choose the unique prime element   $\pi_i$ in $\frak{p}_i$ satisfying $ \pi_i \equiv 1 \; \mbox{mod}\; (3\sqrt{-3}) \; (i = 1, 2).$\\
We set
$$ K_i := k(\sqrt[3]{\pi_i}) = \frak{K}_{\{ \frak{p}_i \}} \;\;\;\; (i = 1, 2).$$
In the following, we assume that
$$ \left(  \frac{\pi_1}{\pi_2} \right)_3 = \left(  \frac{\pi_2}{\pi_1} \right)_3 = 1. \leqno{(5.5)}$$
Note that (5.5)  is  equivalent to that $\frak{p}_1$ (resp. $\frak{p}_2$) splits in $K_2/k$ (resp. $K_1/k$). Let $\frak{P}$ be a fixed prime of $K_1$ lying over $\frak{p}_2$. 
\\
\\
{\bf Proposition 5.6.} {\it  There is an algebraic integer $\alpha$ in $K_1$ which satisfies the following properties.}\\
(1) {\em ${\rm N}_{K_1/k}(\alpha) =  \pi_2 \beta^3$ with some $\beta \in k$.} \\
(2) {\em The principal ideal $(\alpha)$ has the decomposition of the form}
$$ 
(\alpha) = \frak{P}^a \frak{B}^b,  \; \; (a,3) = 1,\;\;  (\frak{B}, 3) = 1, \;\;  b \equiv 0 \; \mbox{mod} \; 3.
$$
\\
{\em Proof.}  (1) Write $\frak{P}^{h_{K_1}} = (\alpha')$ for some $\alpha' \in {\cal O}_{K_1}$. Since $\frak{p}_2$ is completely decomposed in $K_1/k$, we have
$ {\rm N}_{K_1/k}((\alpha')) = {\rm N}_{K_1/k} \frak{P}^{h_{K_1}} = (\pi_2^{h_{K_1}})$
and so  ${\rm N}_{K_1/k}(\alpha') = \varepsilon \pi_2^{h_{K_1}}$ for some $\varepsilon \in {\cal O}_k^{\times}$. By Proposition 5.4, there is $\frak{u} \in {\cal O}_{K_1}^{\times}$ such that ${\rm N}_{K_1/k}(\frak{u}) = \varepsilon$. Letting $\alpha := \frak{u}^{-1} \alpha'$, we have ${\rm N}_{K_1/k}(\alpha) = \pi_2^{h_{K_1}}$. Since $h_{K_1}$ is prime to $3$ by Proposition 5.3, $\alpha$ or $\alpha^2$ satisfies the desired condition.\\
(2) This follows immediately from (1). $\;\; \Box$
\\
\\
Let $\alpha$ be  an element of ${\cal O}_{K_1}$ satisfying the  properties (1), (2) of  Proposition 5.6.
Let $\tau$ be the element of Gal($K_1/k$) defined by 
$$
\tau( \sqrt[3]{\pi_1})  := \zeta_3 \sqrt[3]{\pi_1}.
$$
Then we have Gal($K_1/k$) = $\langle \tau  \; | \; \tau^3 = 1\rangle$.
We set
$$ \left\{ 
\begin{array}{l} \alpha^{(1)} := \alpha,\\
\alpha^{(2)} := \tau(\alpha), \\
\alpha^{(3)} := \tau^2(\alpha),\\
\end{array}\right.
\;\; \;\; \;\; \;\; 
 \left\{
\begin{array}{l} 
\frak{P}^{(1)} := \frak{P}, \\
\frak{P}^{(2)} := \tau(\frak{P}), \\
\frak{P}^{(3)} := \tau^2(\frak{P}),\\ 
\end{array}\right.
$$
where $\frak{P}^{(1)}, \frak{P}^{(2)}$ and $\frak{P}^{(3)}$ are distinct all prime ideals of ${\cal O}_{K_1}$ lying over $\frak{p}_2$.
By Proposition 5.6, we  see easily the following.
\\
\\
{\bf Theorem  5.7.}   {\it Let $\theta := \zeta_3^{c}(\alpha^{(1)})^2 \alpha^{(2)} \in {\cal O}_{K_1} \; (c = 0,1,2)$. Then $\theta$ satisfies the following properties}:\\
(1) ${\rm N}_{K_1/k}(\theta) =  \pi_2^3 w^3$ {\em for some}  $w \in k$. Writing $\theta = x + y \sqrt[3]{\pi_1} + z (\sqrt[3]{\pi_1})^2$ $(x,y,z \in k)$, it is given by
$$ x^3 + \pi_1 y^3 + \pi_1^2 z^3 - 3 \pi_1 xyz = \pi_2^3 w^3.  $$
(2) $(\theta) = (\frak{P}^{(1)})^{2e}(\frak{P}^{(2)})^e \frak{A}^{a},  \; \; (e,3) = 1,  \; \;  (\frak{A}, 3) = 1, \;\; a \equiv 0 \; \mbox{mod}\; 3.$
\\
\\
In the following, we assume that  $\theta$ in Theorem 5.7 satisfies the following condition: \\
\\
(A) $\;\;$ There is $\eta \in {\cal O}_{K_1}$ such that $\eta^3 \equiv \theta \mod (3\sqrt{-3})$. \\ 
\\
A sufficient condition for (A) to hold is given as follows. Set $\frak{U}_{K_1} := ({\cal O}_{K_1}/(3\sqrt{-3}))^{\times}$ and let $\frak{U}_{K_1}(3)$ denote the $3$-Sylow subgroup of $\frak{U}_{K_1}$. By replacing $\alpha$ by $\alpha^{b^2}$, we may assume that $\alpha$ in Proposition 5.6 satisfies $\overline{\alpha} \in \frak{U}_{K_1}(3)$, where $b$ (= 8 or 26) is the non-$3$ part of $\# \frak{U}_{K_1}$. 
\\
\\
{\bf Proposition 5.8.} (1) {\em The group $\frak{U}_{K_1}(3)$ is given by}
$$ \frak{U}_{K_1}(3) = \langle \overline{a_1} \rangle \times \langle \overline{a_2} \rangle \times  \langle \overline{a_3} \rangle \times  \langle \overline{a_4} \rangle \times  \langle \overline{a_5} \rangle \times 
\langle \overline{a_6} \rangle, \;\; \langle \overline{a_i} \rangle \simeq \mathbb{Z}/3\mathbb{Z},$$
{\em where  we set $\overline{a_i} := a_i$ mod $(3\sqrt{-3}), \lambda := \frac{1}{\sqrt{-3}}(\sqrt[3]{\pi_1} -1)$ and}
$$ \left\{ 
\begin{array}{l} 
a_1 := 4, \\
a_2 := \zeta_3,\\
 a_3 := \sqrt[3]{\pi_1}, \\
a_4 := 1 + \sqrt{-3} \sqrt[3]{\pi_1} = 1 + \sqrt{-3} + 3\lambda, \\
a_5 := -1 + 3\sqrt[3]{\pi_1} - (\sqrt[3]{\pi_1})^2 = 1 + \sqrt{-3} \lambda + 3 \lambda^2, \\
a_6 := \frac{1}{\sqrt{-3}} ( 1 + (\sqrt{-3} -2) \sqrt[3]{\pi_1} + (\sqrt[3]{\pi_1})^2) = 1 + \sqrt{-3} \lambda + \sqrt{-3} \lambda^2
\end{array}\right.$$
(2) {\em  Assume that $\alpha$ in Proposition 5.6 satisfies $\overline{\alpha} \in \frak{U}_{K_1}(3)$. Then  $\overline{\alpha} \in \langle \overline{a_1} \rangle \times \langle \overline{a_2} \rangle \times  \langle \overline{a_3} \rangle \times  \langle \overline{a_4} \rangle$. If we have $\overline{\alpha} \in \langle \overline{a_1} \rangle \times \langle \overline{a_2} \rangle \times  \langle \overline{a_3} \rangle$, then 
 $\theta$ in Theorem 5.7 satisfies the condition ${\rm (A)}$ for some $c$.}
\\
\\
{\em Proof.} (1) Since $(\sqrt{-3})$ is unramified in $K_1/k$, we can write $(\sqrt{-3}) = \frak{q}_1 \cdots \frak{q}_r$ ($r=1$ or $3$) with prime ideals $\frak{q}_i$'s of ${\cal O}_{K_1}$. Then we have
$$
\begin{array}{cl}
\frak{U}_{K_1} &\simeq ({\cal O}_{K_1}/\frak{q}_1^3 \cdots \frak{q}_r^3)^{\times} \\
  &\simeq ({\cal O}_{K_1}/\frak{q}_1^3)^{\times} \times \cdots \times ({\cal O}_{K_1}/\frak{q}_r^3)^{\times}.
\end{array}
$$
We note that $\#\frak{U}_{K_1} = \prod_{i=1}^r {\rm N}\frak{q}_i^2 ({\rm N}\frak{q}_i  - 1) = {\rm N}( (\sqrt{-3}{\cal O}_{K_1})^2 ) \prod^{r}_{i=1}({\rm N}\frak{q}_i - 1) = 729b$, where $b =  \prod^{r}_{i=1}({\rm N}\frak{q}_i - 1) \in \{8, 26\}$, and that  $a_i^3 \equiv 1$ mod $(3\sqrt{-3})$. Define $\tau \in$ Gal($K_1/k$), $\tau(\sqrt[3]{\pi_1}) = \zeta_3 \sqrt[3]{\pi_1}$. Then we have 
$$
\tau(a_3) \equiv a_2  \cdot a_3, \; \; \;  \tau(a_4) \equiv a_1 \cdot a_4 \;\; \mbox{mod} \; (3\sqrt{-3}) \\
$$
and hence $\overline{a_3} \not\in \langle \; \overline{a_1}, \overline{a_2} \;\rangle$ and $\overline{a_4} \not\in \prod_{i=1}^3 \langle \; \overline{a_i} \;\rangle$. We note that ${\rm N}_{K_1/k}(a_i) \equiv 1$ mod $(3\sqrt{-3})$ ($1 \le i \le 4$) and
$$
\begin{array}{cl}
{\rm N}_{K_1/k}(a_5) &= -1 + 27 \pi -\pi^2 - 3 \cdot (-1) \cdot 3 \cdot (-1) \cdot \pi \\
         &\equiv -1 -\pi^2 \equiv 7 \mod (3\sqrt{-3}) \\
{\rm N}_{K_1/k}(a_6) &= \frac{\sqrt{-3}}{9}(1 + (\sqrt{-3}-2)^3  \pi_1 + \pi_1^2 - 3 \cdot 1 \cdot (\sqrt{-3}-2) \cdot 1 \cdot \pi) \\
         &= \frac{\sqrt{-3}}{9}(1 + (16 + 6\sqrt{-3}) \pi_1 + \pi_1^2) \\
         &=  \frac{\sqrt{-3}}{9}( (1- \pi_1)^2 + (18 + 6\sqrt{-3}) \pi_1) \\
         &=  \frac{\sqrt{-3}}{9}(1- \pi_1)^2 + (2\sqrt{-3} - 2) \pi_1  \\
         &\equiv (2\sqrt{-3} - 2) \pi_1 \equiv 7 + 2\sqrt{-3} \mod (3\sqrt{-3})
\end{array}
$$
and hence $\overline{a_5} \not\in  \prod_{i=1}^4 \langle \; \overline{a_i}  \; \rangle$, $ \overline{a_6} \not\in  \prod_{i=1}^5 \langle \; \overline{a_i} \; \rangle$.
Since the order of the group $\prod_{i=1}^6 \langle \overline{a_i} \rangle$ is $3^6 = 729$,  it must be $\frak{U}_{K_1}(3)$.   \\
(2) Since ${\rm N}_{K_1/k}(\alpha) \equiv 1$ mod $(3\sqrt{-3})$, by (1), we have $\overline{\alpha} \in \prod_{i=1}^4 \langle \overline{a_i} \rangle$. Suppose $\overline{\alpha} \in \prod_{i=1}^3 \langle \overline{a_i} \rangle$. Then we can write $\overline{\alpha} = \overline{a_1}^{b_1} \cdot \overline{a_2} ^{b_2} \cdot \overline{a_3}^{b_3}$ and $\overline{\tau(\alpha)} = \overline{a_1}^{b_1} \cdot \overline{a_2} ^{b_2 + b_3} \cdot \overline{a_3}^{b_3}$ by (1).  Therefore  we have $\theta \equiv 1$ mod $(3\sqrt{-3})$ with $\theta = \zeta_3^{2b_3} (\alpha^{(1)})^2 \alpha^{(2)}$.  $ \;\; \Box$
\\
\\
{\bf Corollary 5.9.} {\em Assume that $\frak{p}_1=(p_1)$ and $\frak{p}_2=(p_2)$ with prime numbers $p_1$ and $p_2$. 
Then (5.5) is satisfied, and there exists $\alpha$ as in Proposition 5.6 such that $\theta$ in Theorem 5.7 satisfies condition ${\rm (A)}$ for some $c$. }
\\
\\
{\em Proof.} 
Since $\frak{p}_1=(p_1)$ with prime number $p_1$, $K_1/\mathbb Q$ is a Galois extension. 
Note that $p_i$ is inert in $k/\mathbb Q$, and that $\pi_i=-p_i$. 
Since $K_1/\mathbb Q$ is not cyclic, $\frak{p}_2=(p_2)$ splits in $K_1/k$. 
Similarly, $\frak{p}_1$ splits in $K_2/k$. Hence (5.5) is satisfied. 
Put $F := \mathbb Q(\sqrt[3]{p_1})$. 
In the proof of Proposition 5.6, we may choose a prime $\frak{P}$ of $K_1$ such that $F$ is the decomposition field of $\frak{P}$ in $K_1/\mathbb Q$. 
Then $\frak{P}=\wp {\cal O}_{K_1}$ with a prime ideal $\wp$ of ${\cal O}_F$, 
and a principal ideal $\frak{P}^{h_{K_1}} = \wp^{h_{K_1}} {\cal O}_{K_1} = (\alpha')$ has a generator $\alpha' = x+y\sqrt[3]{p_1}+z(\sqrt[3]{p_1})^2 \in F$ with $x,y,z \in \mathbb Q$. 
Since ${\rm N}_{K_1/k}(\alpha') = x^3+p_1y^3+p_1^2z^3-3p_1xyz \in \mathbb Q$, 
we have $\varepsilon=(-p_2)^{-h_{K_1}}{\rm N}_{K_1/k}(\alpha')=\pm 1=\frak{u}$ in the proof of Proposition 5.6. 
Thus we can take $\alpha \in F$ in Proposition 5.6. 
Moreover, we may assume that $\overline{\alpha} \in \frak{U}_{K_1}(3)$. 
Then $\sigma(\alpha)=\alpha$, where $\sigma$ is the generator of $\mathrm{Gal}(K_1/\mathbb Q(\sqrt[3]{p_1}))$. 
By Proposition 5.8 (2), we can write 
$\overline{\alpha} = \overline{a_1}^{b_1} \cdot \overline{a_2} ^{b_2} \cdot \overline{a_3}^{b_3} \cdot \overline{a_4}^{b_4}$. 
We can easily see that
\[
a_1 \cdot a_4^2 \equiv 1-\sqrt{-3}\sqrt[3]{-p_1} = \sigma(a_4) \; \mbox{mod} \; (3\sqrt{-3}) , 
\]
where we note that $\sqrt[3]{-p_1} -1 = \sqrt{-3}\lambda \equiv 0 \; \mbox{mod}\; (\sqrt{-3})$. 
Since
\[
\overline{a_1}^{b_1} \cdot \overline{a_2} ^{b_2} \cdot \overline{a_3}^{b_3} \cdot \overline{a_4}^{b_4} 
= \overline{\alpha} = \sigma(\overline{\alpha}) = 
\overline{a_1}^{b_1+b_4} \cdot \overline{a_2} ^{2b_2} \cdot \overline{a_3}^{b_3} \cdot \overline{a_4}^{2b_4} ,
\]
we have $b_4 \equiv 0 \; \mbox{mod} \;3$, i.e., $\overline{\alpha} \in \prod_{i=1}^3 \langle \overline{a_i} \rangle$. 
By Proposition 5.8 (2), we obtain the claim. $ \;\; \Box$\\
\\
Let $\theta $ be an element of ${\cal O}_{K_1}$ in Theorem 5.7 satisfying  condition (A).  We then set\\
$$
\frak{K}_{\theta} := k(\sqrt[3]{\pi_{1}},\sqrt[3]{\pi_2},\sqrt[3]{\theta}).
$$
\\
{\bf Remark 5.10.}  The extension $\frak{K}_{\theta}$ over $k$ may be regarded as a cubic analogue of R\'{e}dei's dihedral extension $\frak{R}$ over $\mathbb{Q}$ given by
$$ \left\{ \begin{array}{l} \frak{R} = \mathbb{Q}(\sqrt{p_1}, \sqrt{p_2}, \sqrt{\alpha}),  \; \alpha = x + y \sqrt{p_1} \; \mbox{with}\; x, y, z \in \mathbb{Z} \; \mbox{satisfying} \\
x^2 - p_1y^2 = p_2 z^2, (x,y,z) = 1, y \equiv 0 \; {\rm mod}\;  2, x-y \equiv 1 \; {\rm mod}\; 4. \end{array} \right. $$
Here we may observe that the equation  $ x^3 + \pi_1 y^3 + \pi_1^2 z^3 - 3 \pi_1 xyz = \pi_2^3 w^3$ in (1) of Theorem 5.7 corresponds to the equation $x^2 - p_1y^2 = p_2 z^2$, the property (2) of Theorem 5.7 corresponds to $(x,y,z) = 1$  and the condition (A) corresponds to $y \equiv 0$ mod $2$ and $x-y \equiv 1$ mod $4$. \\
\\
In the following, we will show that $\frak{K}_{\theta}$ enjoys the  properties of an  $H_3(\mathbb{F}_3)$-extension over $k$ associated to $\{ \frak{p}_1, \frak{p}_2 \}$.  \\
\\
{\bf Theorem 5.11.} {\em We assume the condition ${\rm (A)}$. Then the following assertions hold.}\\
(1) {\it We have $\frak{K}_{\theta} = k(\sqrt[3]{\theta^{(1)}},\sqrt[3]{\theta^{(2)}},\sqrt[3]{\theta^{(3)}})$, where $\theta^{(1)} := \theta = \zeta_3^{c} (\alpha^{(1)})^2 \alpha^{(2)}, \theta^{(2)} := \tau(\theta) = \zeta_3^{c} (\alpha^{(2)})^2 \alpha^{(3)},$ and $\theta^{(3)} := \tau^2(\theta) = \zeta_3^{c} (\alpha^{(3)})^2 \alpha^{(1)}.$} \\
(2) {\em The extension  $\frak{K}_{\theta}/k$ is Galois. This extension  is unramified outside $\frak{p}_1$ and $\frak{p}_2$ and the ramification index of each $\frak{p}_i$ is $3$. }\\
(3) {\em The Galois group ${\rm Gal}(\frak{K}_{\theta}/k)$ is isomorphic to $H_3(\mathbb{F}_3)$.}\\
\\
{\em Proof.} (1) $k(\sqrt[3]{\theta^{(1)}},\sqrt[3]{\theta^{(2)}},\sqrt[3]{\theta^{(3)}}) \subset \frak{K}_{\theta}$: Since $\sqrt[3]{\pi_2}, \sqrt[3]{\theta^{(1)}} \in \frak{K}_{\theta}$ and we see 
$$\theta^{(1)}(\theta^{(2)})^2 = \pi_2^2(\alpha^{(2)}\beta)^3\;  \mbox{and} \;  \theta^{(1)}\theta^{(2)} = (\theta^{(3)})^2(\alpha^{(2)}/\alpha^{(3)})^3,$$
we have $\sqrt[3]{\theta^{(2)}}, \sqrt[3]{\theta^{(3)}} \in \frak{K}_{\theta}$ and hence the assertion is proved.\\
$\frak{K}_{\theta} \subset k(\sqrt[3]{\theta^{(1)}},\sqrt[3]{\theta^{(2)}},\sqrt[3]{\theta^{(3)}})$: Writing $\theta = x + y\sqrt[3]{\pi_1} + z(\sqrt[3]{\pi_1})^2,$ we have
$$\zeta_3\theta^{(1)} + \theta^{(2)} + \zeta_3^2 \theta^{(3)} = 3 \zeta_3 y \sqrt[3]{\pi_1}.$$
So, $\sqrt[3]{\pi_1} \in k(\sqrt[3]{\theta^{(1)}},\sqrt[3]{\theta^{(2)}},\sqrt[3]{\theta^{(3)}})$ if $y \neq 0$. If $y = 0$, then $z \neq 0$ and $\zeta_3 \theta^{(1)} + \theta^{(2)} = (1 + \zeta_3)x - z(\sqrt[3]{\pi_1})^2$. Hence $\sqrt[3]{\pi_1} \in k(\sqrt[3]{\theta^{(1)}},\sqrt[3]{\theta^{(2)}},\sqrt[3]{\theta^{(3)}})$. Using $\theta^{(1)}(\theta^{(2)})^2 = \pi_2^2(\alpha^{(2)}\beta)^3$ again, we have $\sqrt[3]{\pi_2} \in k(\sqrt[3]{\theta^{(1)}},\sqrt[3]{\theta^{(2)}},\sqrt[3]{\theta^{(3)}})$.  Thus the assertion is proved.\\
(2) Since $\frak{K}_{\theta}$ is the splitting field of $\Pi^{3}_{i=1}(T^3 - \theta^{(i)}) \in {\cal O}_k[T]$, $\frak{K}_{\theta}$ is a Galois extension over $k$. 
We shall show that only $\frak{p}_1$ and $\frak{p}_2$ are ramified in $\frak{K}_{\theta}/k$ with ramification index $3$.
Since $\sqrt[3]{\theta} \notin K_1$ by Theorem 5.7 (2),  we have  $[K_1(\sqrt[3]{\theta}) : K_1] = 3$.
Let $\xi := \sqrt{-3}(\eta - \sqrt[3]{\theta})/3$.
Then $\xi \in {\cal O}_{K_1(\sqrt[3]{\theta})}$, since $\xi$ satisfies $\xi^3 - \sqrt{-3} \eta \xi^2 - \eta^2 \xi + \sqrt{-3}(\eta^3-\theta)/9 = 0$ and $\sqrt{-3}(\eta^3-\theta)/9 \in {\cal O}_{K_1}$ by $\eta^3 \equiv \theta \pmod{3\sqrt{-3}}$(Assumption (A) ).
The relative discriminant of $\xi$ in $K_1(\sqrt[3]{\theta})/K_1$ is computed as 
$$
d(\xi, K_1(\sqrt[3]{\theta})/K_1) = 
\left|
\begin{array}{ccc}
1 & \xi^{(1)} & (\xi^{(1)})^2 \\
1 & \xi^{(2)} & (\xi^{(2)})^2 \\
1 & \xi^{(3)} & (\xi^{(3)})^2 \\
\end{array}
\right|^2
= -\frac{\theta^2}{27}
\left|
\begin{array}{ccc}
1 & 1 & 1 \\
1 & \zeta_3 & \zeta_3^2 \\
1 & \zeta_3^2 & \zeta_3 \\
\end{array}
\right|^2
= \theta^2,
$$
where $\xi^{(1)} := \xi$, $\xi^{(2)} := \sqrt{-3}(\eta - \zeta_3 \sqrt[3]{\theta})/3$ and $\xi^{(3)} := \sqrt{-3}(\eta - \zeta_3^2 \sqrt[3]{\theta})/3$.
 So, only $\frak{P}^{(1)}$ and $\frak{P}^{(2)}$ are ramified in $K_1(\sqrt[3]{\theta})/K_1$.
Similarly, only $\frak{P}^{(2)}$ and $\frak{P}^{(3)}$ are ramified in $K_1(\sqrt[3]{\theta^{(2)}})/K_1$.
Since $\frak{K}_{\theta} = K_1(\sqrt[3]{\theta^{(1)}})\cdot K_1(\sqrt[3]{\theta^{(2)}})$ and only $\frak{p}_1$ is ramified in $K_1/k$, we conclude that only $\frak{p}_1$ and $\frak{p}_2$ are ramified in $\frak{K}_{\theta}/k$ and their ramification indices are 3.
\\
(3)  First, we show $[\frak{K}_{\theta}: k] = 27$. By (2), $K_1(\sqrt[3]{\theta})/K_1$ is a cyclic extension of degree $3$ where only 
$\frak{P}^{(1)}$ and $\frak{P}^{(2)}$ are ramified. By Theorem 3.5, $K_1(\sqrt[3]{\pi_2})/K_1$ is a cyclic extension of degree $3$ where only $\frak{P}^{(i)}$ ($1 \le i \le 3$) are ramified. So, 
$K_1(\sqrt[3]{\theta}) \cap K_1(\sqrt[3]{\pi_2}) = K_1$. Since $\frak{K}_{\theta} = K_1(\sqrt[3]{\theta})\cdot K_1(\sqrt[3]{\pi_2})$, 
$[\frak{K}_{\theta} : K_1] = [ K_1(\sqrt[3]{\theta}) : K_1][K_1(\sqrt[3]{\pi_2}) : K_1] = 9$. Hence $[\frak{K}_{\theta} : k] = [\frak{K}_{\theta} : K_1] [K_1 : k] = 27$.

By the computer calculation using [GAP], we have the following presentation of the group $H_3(\mathbb{F}_3)$:
$$ \begin{array}{ll} 
H_3(\mathbb{F}_3) & =  \left\langle g_1, g_2, g_3  \, \left| \, 
\begin{array}{lr} g_1^3 = g_2^3 = g_3^3 = 1 & \\
 g_3g_2 g_1 = g_1 g_2, \; g_3 g_1 =g_1 g_3, \; g_3 g_2 =g_2 g_3 & \\
\end{array}
\right\rangle 
\right. \\
& = \left\langle g_1, g_2  \; |\;  g_1^3 = g_2^3 = (g_1g_2^2)^3 = (g_1^2 g_2^2)^3 = 1 \right\rangle,
\end{array} 
$$
where $g_1$, $g_2$ and $g_3$ are words representing the following matrices respectively:
$$
g_1 = \left (
\begin{array}{ccc}
1 & 1 & 0  \\ 
0 & 1 & 0  \\
0 & 0 & 1  \\
\end{array}
\right), \;\;
g_2 = \left (
\begin{array}{ccc}
1 & 0 & 0  \\ 
0 & 1 & 1  \\
0 & 0 & 1  \\
\end{array}
\right), \;\;
g_3 = \left (
\begin{array}{ccc}
1 & 0 & 1  \\ 
0 & 1 & 0  \\
0 & 0 & 1  \\
\end{array}
\right) (= [g_1, g_2]).
$$
On the other hand, we define $\gamma_1, \gamma_2, \gamma_3 \in \mbox{Gal}(\frak{K}_{\theta}/k)$ by 
$$
\left.
\begin{array}{rl}
\gamma_1 : & (\sqrt[3]{\pi_1},\sqrt[3]{\pi_2},\sqrt[3]{\theta^{(1)}},\sqrt[3]{\theta^{(2)}},\sqrt[3]{\theta^{(3)}}) \\
 &\mapsto (\zeta_3 \sqrt[3]{\pi_1},\sqrt[3]{\pi_2}, \sqrt[3]{\theta^{(2)}}, \sqrt[3]{\theta^{(3)}}, \sqrt[3]{\theta^{(1)}}) \\
\gamma_2 : &(\sqrt[3]{\pi_1},\sqrt[3]{\pi_2},\sqrt[3]{\theta^{(1)}},\sqrt[3]{\theta^{(2)}},\sqrt[3]{\theta^{(3)}}) \\
 &\mapsto (\sqrt[3]{\pi_1}, \zeta_3 \sqrt[3]{\pi_2},\sqrt[3]{\theta^{(1)}},\zeta_3^2 \sqrt[3]{\theta^{(2)}},\zeta_3 \sqrt[3]{\theta^{(3)}}) \\
\gamma_3 : & (\sqrt[3]{\pi_1},\sqrt[3]{\pi_2},\sqrt[3]{\theta^{(1)}},\sqrt[3]{\theta^{(2)}},\sqrt[3]{\theta^{(3)}}) \\
 & \mapsto (\sqrt[3]{\pi_1},\sqrt[3]{\pi_2}, \zeta_3 \sqrt[3]{\theta^{(1)}},\zeta_3 \sqrt[3]{\theta^{(2)}},\zeta_3 \sqrt[3]{\theta^{(3)}}).
\end{array}
\right.  \leqno{(5.11.1)}
$$
Then we have 
$$ \gamma_1^3 = \gamma_2^3 = \gamma_3^3  = \mbox{id},   \gamma_3\gamma_2 \gamma_1 = \gamma_1 \gamma_2, \; \gamma_3 \gamma_1 = \gamma_1 \gamma_3  \; \mbox{and}\; \gamma_3 \gamma_2 = \gamma_2 \gamma_3,$$
equivalently,
$$ \gamma_3 = [\gamma_1, \gamma_2], \; \gamma_1^3 = \gamma_2^3 = (\gamma_1 \gamma_2^2)^3 = (\gamma_1^2 \gamma_2^2)^3 = 1. \leqno{(5.11.2)} $$
Thus the correspondence $g_i \mapsto \gamma_i$ $(i = 1,2)$ gives a homomorphism  $\kappa : H_3(\mathbb{F}_3) \rightarrow {\rm Gal}(\frak{K}_{\theta}/k)$. Since we easily see that the fixed subfields of $\frak{K}_{\theta}$ by $\langle \gamma_2 \rangle$ and $\langle \gamma_3^2 \rangle$ are $K_1(\sqrt[3]{\theta})$ and $K_1(\sqrt[3]{\pi_2})$, respectively, ${\rm Gal}(\frak{K}_{\theta}/k)$ is generated by $\gamma_1$ and $\gamma_2$, and so $\kappa$ is surjective. Since $H(\mathbb{F}_3)$ and ${\rm Gal}(\frak{K}_{\theta}/k)$ have the same order $27$, $\kappa$ is an isomorphism. $\;\; \Box$
\\
\\
All subgroups of ${\rm Gal}(\frak{K}_{\theta}/k)$ and the corresponding subfields  of $\frak{K}_{\theta}$ are illustrated as follows.\\
\\

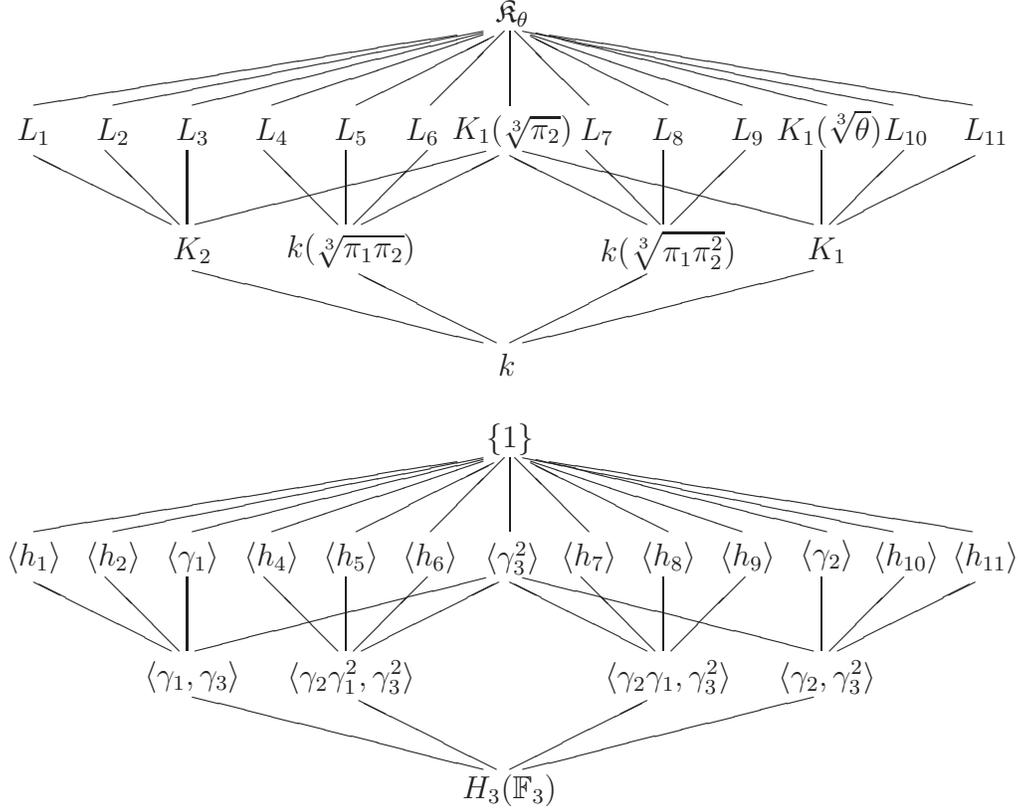
$$
\begin{picture}(100,50)

\put(39,-58){\makebox(20,20){$k$}}
\put(-20,-15){\makebox(20,20){$k(\sqrt[3]{\pi_1\pi_2})$}}
\put(-80,-15){\makebox(20,20){$K_2$}}
\put(100,-15){\makebox(20,20){$k(\sqrt[3]{\pi_1\pi_2^2})$}}
\put(160,-15){\makebox(20,20){$K_1$}}

\put(41,30){\makebox(20,20){$K_1(\sqrt[3]{\pi_2}) $}}

\put(-140,30){\makebox(20,20){$L_1$}}
\put(-110,30){\makebox(20,20){$L_2$}}
\put(-80,30){\makebox(20,20){$L_3$}}
\put(-50,30){\makebox(20,20){$L_4$}}
\put(-20,30){\makebox(20,20){$L_5$}}
\put(7,30){\makebox(20,20){$L_6$}}

\put(73,30){\makebox(20,20){$L_7$}}
\put(100,30){\makebox(20,20){$L_8$}}
\put(130,30){\makebox(20,20){$L_9$}}
\put(161,31){\makebox(20,20){$K_1(\sqrt[3]{\theta})$}}
\put(190,30){\makebox(20,20){$L_{10}$}}
\put(220,30){\makebox(20,20){$L_{11}$}}

\put(41,75){\makebox(20,20){$\frak{K}_{\theta}$}}

\put(40,-40){\line(-4,1){110}}
\put(45,-40){\line(-2,1){52}}
\put(50,-40){\line(2,1){52}}
\put(55,-40){\line(4,1){110}}

\put(-12,5){\line(0,1){28}}
\put(-72,5){\line(0,1){28}}
\put(108,5){\line(0,1){28}}
\put(168,5){\line(0,1){28}}
\put(50,50){\line(0,1){28}}

\put(-15,5){\line(-1,1){28}}
\put(-9,5){\line(1,1){28}}
\put(-6,5){\line(2,1){52}}

\put(111,5){\line(1,1){28}}
\put(107,5){\line(-1,1){28}}
\put(102,5){\line(-2,1){52}}

\put(-69,5){\line(4,1){110}}
\put(-75,5){\line(-1,1){28}}
\put(-78,5){\line(-2,1){52}}

\put(165,5){\line(-4,1){110}}
\put(171,5){\line(1,1){28}}
\put(174,5){\line(2,1){52}}

\put(80,50){\line(-1,1){28}}
\put(110,50){\line(-2,1){52}}
\put(138,50){\line(-3,1){80}}
\put(170,50){\line(-4,1){110}}
\put(195,50){\line(-5,1){140}}
\put(225,50){\line(-6,1){160}}

\put(20,50){\line(1,1){28}}
\put(-8,50){\line(2,1){52}}
\put(-40,50){\line(3,1){80}}
\put(-70,50){\line(4,1){110}}
\put(-100,50){\line(5,1){140}}
\put(-130,50){\line(6,1){160}}

\end{picture}
$$

\vspace{3.2cm}
$$
\begin{picture}(100,50)

\put(40,-58){\makebox(20,20){$H_3(\mathbb{F}_3)$}}
\put(-80,-15){\makebox(20,20){$\langle \gamma_1, \gamma_3 \rangle$}}
\put(-20,-15){\makebox(20,20){$\langle \gamma_2\gamma_1^2, \gamma_3^2 \rangle$}}
\put(100,-15){\makebox(20,20){$\langle \gamma_2\gamma_1, \gamma_3^2 \rangle$}}
\put(160,-15){\makebox(20,20){$\langle \gamma_2, \gamma_3^2 \rangle$}}

\put(41,30){\makebox(20,20){$\langle \gamma_3^2 \rangle$}}

\put(-140,30){\makebox(20,20){$\langle h_1 \rangle$}}
\put(-110,30){\makebox(20,20){$\langle h_2 \rangle$}}
\put(-80,30){\makebox(20,20){$\langle \gamma_1 \rangle$}}
\put(-50,30){\makebox(20,20){$\langle h_4 \rangle$}}
\put(-20,30){\makebox(20,20){$\langle h_5 \rangle$}}
\put(10,30){\makebox(20,20){$\langle h_6 \rangle$}}

\put(70,30){\makebox(20,20){$\langle h_7 \rangle$}}
\put(100,30){\makebox(20,20){$\langle h_8 \rangle$}}
\put(130,30){\makebox(20,20){$\langle h_9 \rangle$}}
\put(160,31){\makebox(20,20){$\langle \gamma_2 \rangle$}}
\put(190,30){\makebox(20,20){$\langle h_{10} \rangle$}}
\put(220,30){\makebox(20,20){$\langle h_{11} \rangle$}}

\put(40,75){\makebox(20,20){$ \{ 1 \}$}}

\put(40,-40){\line(-4,1){110}}
\put(45,-40){\line(-2,1){52}}
\put(50,-40){\line(2,1){52}}
\put(55,-40){\line(4,1){110}}

\put(-12,5){\line(0,1){28}}
\put(-72,5){\line(0,1){28}}
\put(108,5){\line(0,1){28}}
\put(168,5){\line(0,1){28}}
\put(50,50){\line(0,1){28}}

\put(-15,5){\line(-1,1){28}}
\put(-9,5){\line(1,1){28}}
\put(-6,5){\line(2,1){52}}

\put(111,5){\line(1,1){28}}
\put(105,5){\line(-1,1){28}}
\put(102,5){\line(-2,1){52}}

\put(-69,5){\line(4,1){110}}
\put(-75,5){\line(-1,1){28}}
\put(-78,5){\line(-2,1){52}}

\put(165,5){\line(-4,1){110}}
\put(171,5){\line(1,1){28}}
\put(174,5){\line(2,1){52}}

\put(80,50){\line(-1,1){28}}
\put(110,50){\line(-2,1){52}}
\put(138,50){\line(-3,1){80}}
\put(170,50){\line(-4,1){110}}
\put(195,50){\line(-5,1){140}}
\put(225,50){\line(-6,1){160}}

\put(20,50){\line(1,1){28}}
\put(-8,50){\line(2,1){52}}
\put(-40,50){\line(3,1){80}}
\put(-70,50){\line(4,1){110}}
\put(-100,50){\line(5,1){140}}
\put(-130,50){\line(6,1){160}}
\end{picture}
$$
\vspace{1.7cm} \\
where $h_i$ are words representing the following elements respectively:
$$
\begin{array}{c}
h_1 = \gamma_1 \gamma_3^2, \; \; h_2 = \gamma_1 \gamma_3, \; \; h_4 = \gamma_2 \gamma_1^2 , \; \; h_5 = \gamma_2 \gamma_1^2  \gamma_3^2, \; \; h_6 = \gamma_2 \gamma_1^2  \gamma_3, \\
h_7 = \gamma_2 \gamma_1, \; \; h_8 = \gamma_2 \gamma_1 \gamma_3^2, \; \;h_9 = \gamma_2 \gamma_1 \gamma_3, \; \;h_{10} = \gamma_2 \gamma_3^2, \; \;h_{11} = \gamma_2 \gamma_3.
\end{array}
$$
\\
{\bf Corollary 5.12.} {\em Assume condition ${\rm (A)}$. Then we have}
$$ \frak{K}_{\{\frak{p}_1, \frak{p}_2\}} = \frak{K}_{\theta}.$$
{\em In particular, $\frak{K}_{\theta}$ is independent of the choice of $\theta$.}\\
\\
{\em Proof.} This follows from Theorem 4.1 and Theorem 5.11. $\;\; \Box$\\
\\
Corollary 5.9 and Corollary 5.12 yield a constructive proof of Proposition 4.2.\\

\begin{center}
{\bf 6. Triple cubic residue symbols}
\end{center}

In this section, we  introduce the triple cubic residue symbol in terms of the mod $3$ triple Milnor invariant. Then we interpret it arithmetically by using the concrete construction $\frak{K}_{\theta}$ of $\frak{K}_{\{ \frak{p}_1,\frak{p}_2\}}$ in Section 5. It may be regarded as a cubic generalization of R\'{e}dei's triple symbol. \\

Let $\frak{p}_1, \frak{p}_2$ and $ \frak{p}_3$ be distinct primes in  $S_k^{1\, \mbox{{\scriptsize mod}}\, 9}$ and let $\pi_i$ be the unique prime element  in $\frak{p}_i$ satisfying $ \pi_i \equiv 1 \; \mbox{mod} \;  (3\sqrt{-3}) \;  (i = 1, 2, 3)$ as in Lemma 1.1.
We assume that
$$ \left( \frac{\pi_i}{\pi_j} \right)_3 = 1 \;\;\;\; (1\leq i\neq j \leq 3).    \leqno{(6.1)}$$
Moreover, we assume that there exists a R\'{e}dei type $H_3(\mathbb{F}_3)$-extension $\frak{K}_{\{\frak{p}_1, \frak{p}_2\}}$ of $k$ for $\{\frak{p}_1, \frak{p}_2\}$.
By Theorem 4.4, the mod $3$ triple Milnor invariant $\mu_3(123)$ is well-defined.\\
\\
{\bf Definition 6.2.} We define the {\em triple cubic residue symbol} $[\frak{p}_1, \frak{p}_2, \frak{p}_3]_3$ of $\frak{p}_1, \frak{p}_2$ and $ \frak{p}_3$ by
$$ [\frak{p}_1, \frak{p}_2, \frak{p}_3]_3  = \zeta_3^{\mu_3(123)}.$$
\\

We shall describe $[\frak{p}_1, \frak{p}_2, \frak{p}_3]_3$ arithmetically, by using the concrete construction given in Section 5 of the R\'{e}dei type $H_3(\mathbb{F}_3)$-extension $\frak{K}_{\{\frak{p}_1, \frak{p}_2 \}}$ of $k$ for $\{ \frak{p}_1, \frak{p}_2 \}$. Let $\theta \in {\cal O}_{K_1}$ be as in Theorem 5.7 and we assume the condition (A), so that we have
$$\frak{K}_{\{\frak{p}_1, \frak{p}_2 \}} = \frak{K}_{\theta} := k(\sqrt[3]{\pi_{1}},\sqrt[3]{\pi_2},\sqrt[3]{\theta}).$$
Let $\widetilde{\frak{P}}_3$ be a prime of $\frak{K}_3$  lying over $\frak{p}_3$.
Since $\widetilde{\frak{P}}_3$ is unramified in $\frak{K}_{ \{ \frak{p}_1, \frak{p}_2 \}}/k$, we have the Frobenius automorphism $\left(\frac{ \frak{K}_{\{\frak{p}_1, \frak{p}_2 \}}/k }{ \widetilde{\frak{P}}_3 } \right) \in \mbox{Gal}(\frak{K}_{\{\frak{p}_1, \frak{p}_2 \}}/k).$
Since $\frak{p}_3$ splits completely in $k(\sqrt[3]{\pi_1},\sqrt[3]{\pi_2})$ by the assumption (6.1),  we have $ \left(\frac{ \frak{K}_{\{\frak{p}_1, \frak{p}_2 \}}/k }{ \widetilde{\frak{P}}_3 } \right) \in \langle \gamma_3 \rangle.$ Since $\gamma_3 $ is in the center of  ${\rm Gal}(\frak{K}_{\{\frak{p}_1, \frak{p}_2 \}}/k)$,  the Frobenius automorphism $\left(\frac{\frak{K}_{\{\frak{p}_1, \frak{p}_2 \}}/k }{ \widetilde{\frak{P}}_3 } \right)$ is independent of the choice of $\widetilde{\frak{P}}_3$ lying over $\frak{p}_3$ and so we denote it by  
$\left(\frac{\frak{K}_{\{\frak{p}_1, \frak{p}_2 \}}/k }{ \frak{p}_3 } \right).$ 
\\
\\
{\bf Theorem 6.3.} {\em Notations being as above, we have}
$$
\begin{array}{ll}
[\frak{p}_1, \frak{p}_2, \frak{p}_3]_3  & = \displaystyle{ \frac{  \left(  \frac{\frak{K}_{\{\frak{p}_1, \frak{p}_2 \}}/k}{\frak{p}_3}   \right)(\sqrt[3]{\theta} )}{\sqrt[3]{\theta}} }\\
 & = \left\{
\begin{array}{rl}
1 &\quad \left( \frac{\frak{K}_{\{\frak{p}_1, \frak{p}_2 \}}/k}{\frak{p}_3} \right) = \mbox{id}_{\frak{R}_{\{ \frak{p}_1,\frak{p}_2\}}}, \\
\zeta_3 &\quad \left( \frac{\frak{K}_{\{\frak{p}_1, \frak{p}_2 \}}/k}{\frak{p}_3} \right) = \gamma_3 = [\gamma_1, \gamma_2], \\
\zeta_3^2 &\quad \left( \frac{\frak{K}_{\{\frak{p}_1, \frak{p}_2 \}}/k}{\frak{p}_3} \right) = \gamma_3^2 = [\gamma_2, \gamma_1],
\end{array}
\right. 
\end{array}
$$
{\em where $\gamma_1, \gamma_2 , \gamma_3 \in {\rm Gal}(\frak{K}_{ \{ \frak{p}_1, \frak{p}_2 \}}/k)$ are as in (5.11.1).}\\
\\
{\it Proof.} The latter equality  follows from the definition (5.11.1) of $\gamma_i$'s.  We shall prove the former equality.  
Let us  choose a prime $\frak{p}_0  \in S_k^{4, 7\, \mbox{{\scriptsize mod}}\, 9}$ and set $S := \{ \frak{p}_0, \frak{p}_1, \frak{p}_2, \frak{p}_3\}$. Then, by Theorem 1.10, the Galois group $\frak{G}_S$ of the maximal pro-$3$ extension $k_S$ over $k$ unramified outside $S$ has the following presentation:
$$\frak{G}_{S} = \langle x_1, x_2, x_3  \mid  x_{1}^{{\rm N}\frak{p}_1- 1}[x_1 ,y_1] =   x_{2}^{{\rm N}\frak{p}_2- 1}[x_2 ,y_2] = x_{3}^{{\rm N}\frak{p}_3- 1}[x_3 ,y_3] = 1 \rangle,$$
where $x_i$ is the letter representing a monodromy $\tau_i$ over $\frak{p}_i$ in $k_S/k$ and $y_i$ is the free pro-$3$ word of $x_1, x_2, x_3$ representing a Frobenius automorphism $\sigma_i$ over $\frak{p}_i$ in $k_S/k$. 
Let $\frak{F}_3$ be the free pro-3 group on $x_1, x_2, x_3$ and let $\psi : \frak{F}_3  \rightarrow \frak{G}_S $ be the natural homomorphism. 
Since $ \frak{K}_{\{ \frak{p}_1, \frak{p}_2 \}} \subset k_S$, we have the quotient homomorphism $\Phi_S : \frak{G}_{S}  \rightarrow \mbox{Gal}(\frak{K}_{\{ \frak{p}_1, \frak{p}_2 \}}/k) $.
Let $\varphi := \Phi_S \circ \psi  : \frak{F}_3  \rightarrow \mbox{Gal}(\frak{K}_{ \{\frak{p}_1, \frak{p}_2 \} }/k)$. Since we may assume $\tau_i|_{\frak{K}_{\{ \frak{p}_1, \frak{p}_2 \}}} = \gamma_i \, (i=1, 2)$ by Theorem 2.4, we have
$$
\varphi(x_1) = \gamma_1, \quad \varphi(x_2) = \gamma_2, \quad \varphi(x_3) = 1, 
$$
and we have 
$$  \varphi(y_3) = \left( \frac{\frak{K}_{\{ \frak{p}_1, \frak{p}_2 \}} /k}{\frak{p}_3} \right).$$
By (5.11.2),  the relations among $\gamma_1$ and $\gamma_2$ are equivalent to the following relations:
$$ \begin{array}{c} 
\varphi(x_1)^3 = \varphi(x_2)^3 = 1,\quad \varphi(x_3) = 1,\\
 \varphi(x_1 x_2^2)^3 = \varphi(x_1^2 x_2^2)^3 = 1,
\end{array}
$$
and so Ker($\varphi$) is generated as a normal subgroup of $\frak{F}_3$ by 
$$ x_1^3, x_2^3, x_3,   (x_1 x_2^2)^3 \; {\rm and}\; (x_1^2 x_2^2)^3. 
$$
Let  $\Theta_{3} : \frak{F}_3 \rightarrow \mathbb{F}_3\langle \langle X_1, X_2, X_3 \rangle \rangle^{\times}$ be the Magnus embedding of $\frak{F}_3$ over $\mathbb{F}_3$. Then we have
$$ \begin{array}{l}
\quad \Theta_{3}( (x_1)^3 ) = (1 + X_1)^3  = 1 + X_1^3,\\
\quad \Theta_{3}( (x_2)^3 ) = (1 + X_2)^2  = 1 + X_2^3, \\
\quad \Theta_{3}( (x_3) ) = 1 + X_3,  \\
\quad \Theta_{3}( (x_1 x_2^2)^3 ) = ( (1 + X_1) (1 + X_2)^2 )^3   \equiv 1 \; {\rm mod}\;  {\rm deg}\;  \ge 3,\\
\quad \Theta_{3}( (x_1^2 x_2^2)^3 ) = ( (1 + X_1)^2 (1 + X_2)^2 )^3 \equiv 1 \; {\rm mod}\;  {\rm deg}\;  \ge 3.\end{array}
$$
Therefore $\mu_3( (1) ; *),  \mu_3( (2) ; *), \mu_3( (12) ; *)$ take their values 0 on Ker($\varphi$). \\
Case $\left( \frac{\frak{K}_{\{ \frak{p}_1, \frak{p}_2 \}} /k}{\frak{p}_3} \right) = \mbox{id}$:    Then $\varphi(y_3) = 1$ and $\mu_3(123) = \mu_3( (12) ; y_3) = 0$ by $y_3 \in \mbox{Ker}(\varphi)$. \\
Case  $\left( \frac{\frak{K}_{\{ \frak{p}_1, \frak{p}_2 \}}/k}{\frak{p}_3} \right) =  [\gamma_1, \gamma_2]$: Then $\varphi(y_3) = [\gamma_1, \gamma_2]  = \varphi([x_1,x_2])$ and so  we can write $y_3 = [x_1,x_2]  f$ for some  $f \in \mbox{Ker}(\varphi)$. Comparing the coefficients of $X_1X_2$ in the equality $\Theta_{3}(y_3) = \Theta_{3}([x_1,x_2])\Theta_{3}(f)$, we have
$$ \begin{array}{ll}
\mu_3(123) & = \mu_3( (12) ; y_3) \\
            & = \mu_3( (12) ; [x_1,x_2]) + \mu_3( (1) ; [x_1,x_2]) \mu_3( (2) ; f) + \mu_3( (12) ; f) \\
            & = \mu_3( (12) ; [x_1,x_2])\\
 & = 1. 
\end{array}
$$
Case  $\left( \frac{\frak{K}_{\{ \frak{p}_1, \frak{p}_2 \}}/k}{\frak{p}_3} \right) =  [\gamma_2, \gamma_1]$: Then $\varphi(y_3) = [\gamma_2, \gamma_1] = \varphi([x_2,x_1])$ and so we can write $y_3 = [x_2,x_1] f'$ for some $f' \in \mbox{Ker}(\varphi)$. Then comparing the coefficients of $X_1X_2$ in the equality $\Theta_{3}(y_3) = \Theta_{3}([x_2,x_1])\Theta_{3}(f')$, we have
$$ \begin{array}{ll}
\mu_3(123) & = \mu_3( (12) ; y_3) \\
            & = \mu_3( (12) ; [x_2,x_1]) + \mu_3( (1) ; [x_2,x_1]) \mu_3( (2) ; f') + \mu_3( (12) ; f') \\
            & = \mu_3( (12) ; [x_2,x_1]) \\
 & = -1 = 2. \;\; \Box
\end{array}
$$
\\
{\bf Example 6.4.} Let $(\pi_1, \pi_2) := (-17, -53)$.
Then we have
$$ \left\{ \begin{array}{l}
\alpha^{(1)} = \alpha  = 8 - 3 \sqrt[3]{17},\\
\alpha^{(2)}  = 8 - 3 \zeta_3 \sqrt[3]{17},\\
\alpha^{(3)} = 8 - 3 \zeta_3^2 \sqrt[3]{17},\\
\end{array}\right.
\; \; \;\;
\left\{ \begin{array}{l}
\theta^{(1)} = \theta  = (\alpha^{(1)})^2 \alpha^{(2)},\\
\theta^{(2)} = (\alpha^{(2)})^2 \alpha^{(3)}, \\
\theta^{(3)} = (\alpha^{(3)})^2 \alpha^{(1)}.\\
\end{array}\right.
$$
Then we easily see that $\theta$ satisfies (1), (2) in Theorem 5.7. Since  $\alpha^{(1)} \equiv \alpha^{(2)}$ mod $(3\sqrt{-3})$, $\alpha^3 \equiv \theta$ mod $(3\sqrt{-3})$
and so the condition (A) is also satisfied.
Hence
$$\frak{K}_{\{ \frak{p}_1,\frak{p}_2\}} = k(\sqrt[3]{\pi_1},\sqrt[3]{\pi_2},\sqrt[3]{\theta}). $$
Suppose $\pi_3 = -71, -89, -107, -179, -197$. Then we have \\
\\
$ [(17), (53), (71)]_3 = \zeta_3^2$, $[(17), (53), (89)]_3 = \zeta_3$, $[(17), (53), (107)]_3 = \zeta_3^2,$\\
 $[(17), (53), (179)]_3 = \zeta_3$, $[(17), (53), (197)]_3 = \zeta_3$. \\
\\
The right hand side of the equality in Theorem 6.3 depends on the choice of $\theta$ and the order of  $\frak{p}_1$ and $\frak{p}_2$. The mod $3$ Milnor invariant $\mu_3(213)$ is well defined as in Section 4 and hence the triple cubic residue symbol  $[\frak{p}_2, \frak{p}_1, \frak{p}_3]_3$ is defined. Then the shuffle relation (2.1.2) yields the following. \\
\\
{\bf Proposition 6.5.} {\em Notations and assumptions being as above, we have}
$$ [\frak{p}_1, \frak{p}_2, \frak{p}_3]_3  = [\frak{p}_2, \frak{p}_1, \frak{p}_3]_3^{-1} .$$
\vspace{.2cm}

\begin{center}
{\bf 7. Triple Massey products}
\end{center}

In this section, we give a cohomological interpretation of the triple cubic residue symbols in terms of the triple Massey product in Galois cohomology.  Our theorem is seen as a generalization of the known relation between the cup product and the cubic residue symbol to the triple case, and also a generalization of the previous result by the second author ([Mo3], [V]) in the case of the rational number field to the Eisenstein number field. It may be regarded as a mod $3$ arithmetic analogue of the corresponding topological result due to Turaev ([Tu]).\\

Let $\frak{G}$ be a pro-$l$ group ($l$ being a prime number) and let $R$ be a commutative ring with identity on which $\frak{G}$ acts trivally.  Let $C^j(\frak{G},R)$ be the $R$-module of inhomogeneous $j$-cochains $(j \geq 0)$ of $\frak{G}$ with coefficients in $R$ and we consider the differential graded algebra $(C^{\bullet}(\frak{G},R), d)$, 
where the product structure on $C^{\bullet}(\frak{G},R) = \bigoplus_{j\geq 0} C^j(\frak{G},R)$ is given by the cup product and the differential $d$ is the coboundary operator. Then we have the cohomology $H^*(\frak{G},R) = H^*(C^{\bullet}(\frak{G},R))$ of the pro-$l$ group $\frak{G}$ with coefficients in $R$. In the following, we consider Massey products  in $H^1(\frak{G},R)$. For the sign convention, we follow [Dw].   

 Let  $\chi_1,\dots,\chi_n \in H^1(\frak{G},R)$ $(n\geq 2)$. An $n$-th {\it Massey product}  $\langle \chi_1,\dots,\chi_n \rangle$ is said to be {\it defined} if there is an array 
$$ \Omega = \{\omega_{ij} \in C^1(\frak{G},R) \; | \; 1 \leq i < j \leq n+1, (i,j) \neq (1,n+1) \}$$
such that
$$\left\{ 
\begin{array}{l}
[\omega_{i,i+1}] = \chi_i \;\; (1\leq i \leq n),\\
 \displaystyle{d\omega_{ij} = \sum_{a=i+1}^{j-1} \omega_{ia}\cup \omega_{aj}} \;\; (j \neq i+1).

\end{array}
\right.
$$
Such an array $\Omega$ is called a {\it defining system} for $\langle \chi_1,\dots, \chi_n \rangle$. The value of  $\langle \chi_1,\dots, \chi_n \rangle$ relative to $\Omega$,  denoted by $\langle \chi_1,\dots, \chi_n \rangle_{\Omega}$, is defined by the cohomology class represented by the $2$-cocycle 
$$  \sum_{a=2}^n \omega_{1a} \cup \omega_{a,n+1}.$$
We define the Massey product $\langle \chi_1,\dots,\chi_n \rangle $ to be the subset of $H^2(\frak{G},R)$ consisting of elements $\langle \chi_1,\dots, \chi_n \rangle_{\Omega}$ for some defining system $\Omega$. By convention, $\langle \chi \rangle = 0$. We recall the following basic fact (cf. [Kr]). \\
\\
{\bf Lemma 7.2.}  {\em We have  $\langle \chi_1, \chi_2 \rangle = \chi_1 \cup \chi_2$. For $n \geq 3$, $ \langle \chi_1,\dots, \chi_n \rangle$ is defined and consists of a single element if $\langle \chi_{j_1},\dots, \chi_{j_a} \rangle = 0$ for all proper subsets $\{j_1,\dots, j_a \}$ $(a \geq 2)$ of $\{ 1,\dots ,n\}$.} (In this case, we denote the single element by $\langle \chi_{1},\dots, \chi_{n} \rangle$.)\\

Next, we recall a relation between Massey products and Magnus coefficients. Suppose that $\frak{G}$ is  a finitely generated pro-$l$ group with a minimal presentation 
$$ 1 \longrightarrow \frak{N} \longrightarrow \frak{F}_N \stackrel{\psi}{\longrightarrow} \frak{G} \longrightarrow 1,$$
where $\frak{F}_N$ is a free pro-$l$ group ongenerators  $x_1, \dots, x_N$ with $N = {\rm dim}_{\mathbb{F}_l} H^1(\frak{G}, \mathbb{F}_l)$. We set $\tau_i := \psi(x_i) \;(1\leq i \leq N)$. We assume that $\psi$ induces the isomorphism $\frak{F}_N/\Phi(\frak{F}_N) \simeq \frak{G}/\Phi(\frak{G})$ so that  $\psi$ induces the isomorphism  $\psi^* : H^1(\frak{G}, \mathbb{F}_l) \simeq H^1(\frak{F}_N, \mathbb{F}_l).$ We let 
$${\rm tg} : H^1(\frak{N}, \mathbb{F}_l)^{\frak{G}} \rightarrow H^2(\frak{G},\mathbb{F}_l)$$
 be the transgression defined as follows. For $a \in H^1(\frak{N}, \mathbb{F}_l)^{\frak{G}}$, choose a 1-cochain $b \in C^1(\frak{F}_N, \mathbb{F}_l)$ such that $b|_{\frak{N}} = a$. Since the value $db(f_1,f_2)$, $f_i \in \frak{F}_N$, depends only on the cosets $f_i$ mod $\frak{N}$, $db$ defines a 2-cocyle $c$ of $G$. Then ${\rm tg}(a)$ is defined by the class of $c$. By the Hochschild-Serre spectral sequence,  ${\rm tg}$ is an isomorphism and so we have the dual isomorphism, called the Hopf isomorphism,
$${\rm tg}^{\vee} : H_2(\frak{G}, \mathbb{F}_l) \stackrel{\sim}{\rightarrow} H_1(\frak{N}, \mathbb{F}_l)_{\frak{G}} = \frak{N}/\frak{N}^l[\frak{N}, \frak{F}_N].  $$
Then we have the following Proposition (cf. [St; Lemma 1.5], [V; Theorem A.3]). The proof goes in the same manner as in [Mo3; Theorem 2.2.2].\\
\\
{\bf Proposition 7.3.}  {\it Notations being as above, let $\chi_1,\dots,\chi_n \in H^1(\frak{G}, \mathbb{F}_l)$ $(n \geq 2)$. Let $ f \in \frak{N}$ and set $ \delta := ({\rm tg}^{\vee})^{-1}(f \; {\rm mod}\;  \frak{N}^l[\frak{N},\frak{F}_N])$. Assume that   all Massey products up to length $n-1$ are trivial. Then $\frak{N} \subset \frak{F}_N^{(n)}$ and we have} 
$$ \langle \chi_1,\dots, \chi_n \rangle(\delta)  = \displaystyle{(-1)^{n+1} \sum_{{\scriptstyle I = (i_1\cdots i_n)} \atop {\scriptstyle 1 \leq i_1,\dots , i_n \leq N}} \chi_1(x_{i_1})\cdots \chi_n(x_{i_n})\mu_l(I; f).}
$$
\vspace{0.02cm}

Let us be back in our arithmetic situation and keep the same notations as in Section 4. So let $k$ be the Eisenstein number field and let $S_0 := \{ \frak{p}_1, \frak{p}_2, \frak{p}_3 \}$ be a set of 3 distinct primes in $S_k^{1\, \mbox{{\scriptsize mod}}\, 9 }$. As in Lemma 1.1, we choose the unique prime element $\pi_i$ of $\frak{p}_i$ such that $\pi_i \equiv 1 \; \mbox{mod}\; (3\sqrt{-3})$ $(1\leq i \leq 3)$. Choose a prime $\frak{p}_0 \in S_k^{4, 7 \, \mbox{{\scriptsize mod}}\, 9 }$ and set $S := S_0 \cup \{ \frak{p}_0\}$. By Theorem 1.10, we have the following minimal presentation of the Galois group $\frak{G}_S$ of maximal pro-$3$ extension over $k$ unramified outside $S$
$$ \begin{array}{ll} \frak{G}_{S}  & = \langle \,  x_{1}, x_2, x_{3} \; | \; x_{1}^{{\rm N}\frak{p}_{1} -1}[x_{1},y_{1}] =  x_{2}^{{\rm N}\frak{p}_{2} -1}[x_{2},y_{2}] =x_{3}^{{\rm N}\frak{p}_{3} -1}[x_{3},y_{3}] = 1 \rangle \\
& = \frak{F}_3/\frak{N}_S.  \end{array} $$
Here $x_i$ denotes a letter representing a monodromy $\tau_i$ over $\frak{p}_i$ in $k_S/k$ $(1\leq i \leq 3)$ and  $\frak{F}_3$ denotes the free pro-$3$ group on $x_1, x_3, x_3$. The pro-$3$ word  $y_i$ represents a Frobenius automorphism $\sigma_i$ over $\frak{p}_i$ in $k_S/k$ and  $\frak{N}_S$ denotes the closed subgroup of $\frak{F}_r$ generated normally by $x_{i}^{{\rm N}\frak{p}_{i} -1}[x_{i},y_{i}]$ for $1\leq i \leq 3$. We set $\delta_i := ({\rm tg}^{\vee})^{-1}(x_i^{{\rm N}\frak{p}_i -1}[x_i,y_i])$, where ${\rm tg}^{\vee} : H_2(\frak{G}_S,\mathbb{F}_3) \stackrel{\sim}{\rightarrow} \frak{N}_S/\frak{N}_S^3[\frak{N}_S,\frak{F}_3]$ is the Hopf isomorphism.  Let $\chi_1, \chi_2, \chi_3 \in H^1(\frak{G}_{S},\mathbb{F}_3)$ be the Kronecker dual to the monodromies $\tau_1, \tau_2, \tau_3$, namely, $\chi_i(\tau_j) = \delta_{i,j}$.  
 \\
 \\
  {\bf Proposition 7.4.} {\em Notations being as above,  we have, for $1\leq i \neq j \leq 3$,}
    $$ \left(  \frac{\pi_i}{\pi_j} \right)_3 = \zeta_3^{ \langle \chi_i, \chi_j \rangle(\delta_j)},\;\; \left(  \frac{\pi_j}{\pi_i} \right)_3 = \zeta_3^{-\langle \chi_i, \chi_j \rangle(\delta_i)}.$$
  \\
 {\em Proof.} By Proposition 7.3 and the proof of (2.4.2), we have
 $$ \begin{array}{ll} \langle \chi_i, \chi_j \rangle(\delta_a)  & = -\mu_3((ij); x_a^{{\rm N}\frak{p}_a-1}[x_a,y_a])\\
                                                                                          &  = \left\{
  \begin{array}{ll}
  \mu_3(ij) & a = j,\\
  -\mu_3(ji) & a = i,\\
  0 & a \neq i, j.
  \end{array} \right.
  \end{array}
$$
The assertion follows from Theorem 3.6. $\;\; \Box$\\

We assume that 
$$ \langle \chi_i, \chi_j \rangle  = 0 \;\; (1\leq i \neq j \leq 3),$$
which  is equivalent to the condition
$$ \left(  \frac{\pi_i}{\pi_j} \right)_3 = 1 \;\; (1\leq i \neq j \leq 3)$$
by Proposition 7.4, and we assume that the mod $3$ Milnor invariants $\mu_3(abc)$ ($\{a,b,c\} = \{1,2,3\}$) are well defined (cf. Theorem 4.4). By the definition of Massey products and Lemma 7.2, there are 1-cochains $\omega_{13}, \omega_{24} \in C^1(\frak{G}_S, \mathbb{F}_3)$ such that
$$ \langle \chi_1, \chi_2 \rangle = d\omega_{13}, \; \langle \chi_2, \chi_3 \rangle = d\omega_{24},$$
and we have the triple Massey product $\langle \chi_1, \chi_2, \chi_3 \rangle$ defined by
$$ \langle \chi_1, \chi_2, \chi_3 \rangle = [ \chi_1 \cup \omega_{24} + \omega_{13}\cup \chi_3 ].$$
\\
{\bf Theorem 7.5.} {\em Assume that}
$$ \left(  \frac{\pi_i}{\pi_j} \right)_3 = 1 \;\; (1\leq i \neq j \leq 3).$$
{\em Then we have }
$$ [\frak{p}_1, \frak{p}_2, \frak{p}_3]_3 = \zeta_3^{-\langle \chi_1, \chi_2, \chi_3 \rangle(\delta_3)},\;\; [\frak{p}_2, \frak{p}_3, \frak{p}_1]_3 = \zeta_3^{\langle \chi_1, \chi_2, \chi_3 \rangle(\delta_1)}.$$
\\
{\em Proof.} By Proposition 7.3 and the proof of (2.4.2), we have
 $$ \begin{array}{ll} \langle \chi_1, \chi_2, \chi_3 \rangle(\delta_a)  & = \mu_3((123); x_a^{{\rm N}\frak{p}_a-1}[x_a,y_a])\\
                                                                                          &  = \left\{
  \begin{array}{ll}
  -\mu_3(123) & a = 3,\\
  \mu_3(231) & a = 1,\\
  0 & a \neq 1, 3.
  \end{array}\right.
  \end{array}
$$\\
 By Theorem 6.4, we obtain the assertion. $\;\; \Box$
 \\

{\small
\begin{flushleft}
{\bf References}\\
{ [A1] } F. Amano, On R\'{e}dei's dihedral extension and triple reciprocity law, Proc. Japan Acad.,  {\bf 90}, Ser. A  (2014), 1--5. \\
{ [A2] } F. Amano, On a certain nilpotent extension  over $\mathbb{Q}$ of degree 64 and the 4-th multiple residue symbol, Tohoku Math. J. {\bf 66} No.4 (2014), 501 --522.\\
{ [A3] } F. Amano, Arithmetic of nilpotent extensions and multiple residue symbols, Thesis, Kyushu University, 2014.\\
{ [B] } B. J. Birch, Cyclotomic fields and Kummer extensions, In:Algebraic Number Theory
(Proc. Instructional Conf., Brighton, 1965), pages 85--93. Thompson, Washington, D.C., 1967.\\
{ [CFL] } K. T. Chen, R. H. Fox and R. C. Lyndon, Free differential calculus. IV. The quotient groups of the lower central series, Ann. of Math. (2) {\bf 68} (1958), 81--95.\\
{[DDMS]} J. D. Dixon, M. P. F. du Sautoy, A. Mann, D. Segal, Analytic pro-$p$  groups, Second edition. Cambridge Studies in Advanced Mathematics, {\bf 61}, Cambridge University Press, Cambridge, 1999.\\
{[Dw]} W. G. Dwyer, Homology, Massey products and maps between groups, J. Pure Appl. Algebra  {\bf 6},  (1975), no. 2, 177--190. \\
{ [F\"{o}] } A. Fr\"{o}hlich, A prime decomposition symbol for certain non Abelian number fields, ActaSci. Math., {\bf 21}  (I960), 229--246.\\
{ [Fu] } Y. Furuta, A prime decomposition symbol for a nonabelian central extension which is abelian over a bicyclic biquadratic field,  Nagoya Math. J.  {\bf 79}  (1980), 79--109.\\
{[GAP]} The GAP Group, GAP --- Groups, Algorithms, and Programming, ver. 4.7.8, 2015. \\
{ [G] } C. F. Gauss, Disquisitiones arithmeticae, Translated into English by A. Arthur, S. J. Clarke, Yale University Press, New Haven, Conn.-London 1966.\\
{ [H] } T. Honda, Pure cubic fields whose class numbers are multiples of three, J. Number Theory {\bf 3} (1971), 7--12. \\
 { [Ih] } Y. Ihara, On Galois representations arising from towers of coverings of ${\bf P}^1 \setminus \{0,1,\infty \}$, Invent. Math.
 {\bf 86} (1986), no. 3, 427--459.\\
 { [Iw1] }   K. Iwasawa, A note on class numbers of algebraic number fields, Abh. Math. Sem. Hamburg, {\bf 20} (1956), 257--258.\\
 { [Iw2] }  K. Iwasawa, On Galois groups of local fields, Trans. Amer. Math. Soc.  {\bf 80}  (1955), 448--469. \\
 { [Ka] } M. Kapranov, Analogies between number fields and 3-manifolds, unpublished note (1996), Max Planck Institute.\\
 { [Ko] }  H. Koch, Galoissche Theorie der $p$-Erweiterungen, Springer-Verlag, Berlin-New York; VEB Deutscher Verlag der Wissenschaften, Berlin, 1970.\\
 { [Kr] } D. Kraines, Massey higher products, Trans. Amer. Math. Soc.  {\bf 124}  (1966),  431--449.\\
 {[L]} F. Lemmermeyer, The ambiguous class number formula revisited, J. Ramanujan Math. Soc.
{\bf 28} (2013), no. 4, 415-421.\\
 { [Mi1] } J. Milnor, Link groups, Ann. of Math. {\bf 59} (1954), 177--195.\\
{ [Mi2] } J. Milnor, Isotopy of links, in Algebraic Geometry and Topology, A symposium in honor of S. Lefschetz ( edited by R.H. Fox, D.C. Spencer and A.W. Tucker), 280--306 Princeton University Press, Princeton, N.J., 1957.\\
{ [Mo1] } M. Morishita, Milnor's link invariants attached to certain Galois groups over {\bf Q}, Proc. Japan Acad. Ser, A {\bf 76} (2000), 18--21.\\
{ [Mo2] } M. Morishita, On certain analogies between knots and primes, J. Reine Angew. Math. {\bf 550} (2002), 141--167.\\
{ [Mo3] } M. Morishita, Milnor invariants and Massey products for prime numbers, Compos. Math., {\bf 140} (2004), 69--83. \\
 { [Mo4] } M. Morishita, Knots and Primes - An introduction to arithmetic topology, Universitext, Springer, London, 2012.\\
 { [Mu] } K. Murasugi, Nilpotent coverings of links and Milnor's invariant, Low-dimensional topology (Chelwood Gate, 1982), 106-142 London Math. Soc. Lecture Note Ser., 95, Cambridge Univ. Press, Cambridge-New York 1985.\\
{[O]} T. Oda, Note on meta-abelian quotients of pro-$l$ free groups, (1985), preprint.\\
{[Pari-gp]} The PARI~Group, PARI/GP version 2.7.4, Univ. Bordeaux, 2015. \\
{[Rd]} L. R\'{e}dei, Ein neues zahlentheoretisches Symbol mit Anwendungen auf die Theorie der quadratischen Zahlk\"{o}rper I, J. Reine Angew. Math., {\bf 180} (1939), 1-43.\\
{ [Rz]} A. Reznikov, Embedded incompressible surfaces and homology of ramified coverings of three-manifolds, Selecta Math. (N.S.)  {\bf 6}  (2000),  no. 1, 1--39.\\
 { [Se]} J.-P. Serre,  Corps locaux,  Publications de l'Universit\'{e} de Nancago, No. VIII. Hermann, Paris, 1968.\\
 { [St] } D. Stein, Massey products in the cohomology of groups with applications to link theory, Trans.
Amer. Math. Soc. {\bf 318}  (1990), 301--325.\\
{ [Ta]} J. Tate,  Global class field theory,  In:Algebraic Number Theory
(Proc. Instructional Conf., Brighton, 1965), pages 162 -- 203. Thompson, Washington, D.C., 1967.\\
{ [Tu] } V. Turaev, The Milnor invariants and Massey products, (Russian) Studies in topology, II. Zap. Nau\v{c}n. Sem. Leningrad. Otdel. Mat. Inst. Steklov. (LOMI)  {\bf 66},  (1976), 189--203, 209--210. \\
{ [V] } D. Vogel, On the Galois group of 2-extensions with restricted ramification, J. Reine Angew. Math. {\bf 581} (2005), 117--150.\\
{[Y]} H. Yokoi, On the class number of a relatively cyclic number field, Nagoya Math. J.  {\bf 29}  1967 31--44.\\
\end{flushleft} 
\ \\
\\
Fumiya Amano \\
2077-1, Jyuni-cho, Ibusuki, Kagoshima, 891-0403, JAPAN \\
e-mail: ca-solitudeam.p@ezweb.ne.jp\\
\\
Yasushi Mizusawa\\
Department of Mathematics, Nagoya Institute of Technology \\
 Gokiso, Showa, Nagoya, Aichi, 66-8555 JAPAN\\
e-mail: mizusawa.yasushi@nitech.ac.jp\\
\\
Masanori Morishita\\
Faculty of Mathematics, Kyushu University \\
744, Motooka, Nishi-ku, Fukuoka, 819-0395, JAPAN \\
e-mail: morisita@math.kyushu-u.ac.jp
}
\end{document}